\documentclass[a4paper,11pt]{article}
        
\usepackage{macro_Sam} 
\author{Samuel TAPIE}
\date{08/01/2010}
\title{Generic metrics, eigenfunctions and riemannian coverings of non compact manifolds}

\begin{document}

\maketitle

\selectlanguage{english}

\begin{abstract}
Let $(M^n,g)$ be a non-compact riemannian $n$-manifold with bounded geometry at order $k\geq\frac{n}{2}$. We show that if the spectrum of the Laplacian $\Delta_g$ starts with $q+1$ discrete eigenvalues $\lambda_0<\lambda_1\leq...\leq \lambda_q$ isolated from the essential spectrum, and if the metric is generic for the $\Cl C^{k+2}$-strong topology, then the $\lambda_j$ are distinct and their associated eigenfunctions are Morse. This generalizes to non-compact manifolds some arguments developped by K. Uhlenbeck. We deduce from this result that if $M^n$ has bounded geometry at order $k\geq\frac{n}{2}$ and has an isolated first eigenvalue for its Laplacian, then for any riemannian covering $p : M'\ra M$, we have $\lambda_0(M) = \sup_D \lambda_0(D)$, where $D\subset M'$ runs over all connected fundamental domains for $p$, and $\lambda_0(D)$ is the bottom of the spectrum of $D$ with Neumann boundary conditions.
\end{abstract}

\tableofcontents{}

\section*{Introduction}

K. Uhlenbeck has shown in 1976 that on a compact riemannian manifold, when the metric is generic the eigenvalues of the Laplacian are simple and the associated eigenfunctions are Morse. We extend this result to non-compact manifolds with bounded geometry. Recall that a riemannian manifold has \emph{bounded geometry at order $k$} if its injectivity radius is positive and the $k$ first covariant derivatives of the Riemann curvature tensor are uniformly bounded. Our main result is the following.

\begin{Theo*}\label{th:PhiM}
Let $(M^n, g_0)$ be a riemannian manifold with bounded geometry at order $k\geq\frac{n}{2}$. Assume the Laplacian $\Delta_{g_0}$ has $q+1$ first eigenvalues $\lambda_0<\lambda_1\leq...\leq\lambda_q$ (repeated according to their multiplicity) isolated from the essential spectrum of $M$. Then there is a neighbourhood $\Cl U(g_0)$ for the $\Cl C^{k+2}$-strong topology such that for any metric in $\Cl U(g_0)$, the geometry is bounded at order $k$, and the Laplacian has $q+1$ first eigenvalues isolated from the essential spectrum. Moreover, there is a generic set $\Cl U_M(g_0)\subset \Cl U(g_0)$ such that for any metric in $g\in\Cl U_M(g_0)$, the $q+1$ first eigenvalues of $\Delta_g$ are simple and the associated eigenfunctions are Morse.
\end{Theo*}

Our first section presents this \ind{$\Cl C^{k+2}$-strong topology}, and uses it as a crucial tool to adapt Uhlenbeck's arguments to non-compact manifolds. In our second section, we apply these generic properties of the first eigenfunction to study some spectral properties of riemannian coverings :

\begin{Theo*}\label{th:SpecDom}
Let $M$ be a riemannian $n$-manifold with bounded geometry at order $k\geq\frac{n}{2}$. Assume the bottom of the spectrum of $\Delta_g$ is an isolated eigenvalue. Then for any riemannian covering $p : M'\ra M$, we have
$$\lambda_0(M) = \sup_D\{\lambda_0(D), D\subset M'\mbox{ fundamental domain for } p\},$$
where we set Neumann conditions on $\bd D$. This supremum is attained if there is a fundamental domain $D_0\subset M$ on which the function $p\circ\phi_0$ satisfies Neumann boundary condition, where $\phi_0$ is the first eigenfunction of $\Delta_g$.
\end{Theo*}

To show this equality we approximate a given metric by generic metrics, and we use the following characterization of the bottom of the spectrum with Neumann conditions : 

\begin{Theo*}
Let $D$ be a riemannian manifold with boundary, then the bottom of the spectrum of the Laplacian with Neumann boundary conditions is given by
$$\lambda_0(D) = \sup \left\{\lambda\in\Bb R : \exists \phi\in\Cl C^\infty(D), \phi>0, \left.\frac{\bd \phi}{\bd\nu}\right|_{\bd D} \equiv 0, \Delta \phi = \lambda \phi\right\}.$$
\end{Theo*}

This means for any real $\lambda$, there exists a $\Cl C^\infty$ positive $\lambda$-harmonic function $\phi$ on $D$ with Neumann boundary conditions \emph{if and only if} $\lambda\leq \lambda_0(D)$. D. Sullivan gave a proof of this result, based on a diffusion process, when the boundary is empty. In our Appendix, we extends his proof to the case of Neumann boundary conditions.

\paragraph{Acknowledgements} A large part of this work was done during my PhD in Grenoble. Many thanks to G. Besson and G. Courtois for precious advice and encouragements. Thanks to D. Piau for crucial help in stochastic calculus, to R. Joly for introducing me to the Whitney topology, and to G. Carron for important corrections. This article was finished while being hosted by the Hausdorff Center for Mathematics, Bonn.

\section{Generic metrics and eigenfunctions of the Laplacian}

This section is dedicated to the proof of our theorem of genericity. The two first paragraphs present some background and tools we will need. The third paragraph contains the proof itself, which is the main technical part of this paper.

\subsection{Basics on the spectrum and Sobolev spaces}\label{ssec:BasicSpec}

Let $(M, g)$ be a non compact riemannian $n$-manifold. We will call \ind{Laplacian} the \ind{Laplace-Beltrami operator}, defined on $\Cl C^2$ functions $f : M\ra \Bb R$ by
$$\Delta_g f = \gdiv\nabla_g f = -\Trace(\nabla_g df),$$
where $\nabla_g$ is the Levi-Civita connexion associated to $g$. We will omit the index $g$ when no confusion shall arise. In the neighbourhood of a point where the manifold is defined by smooth coordinates $(x_i)_{i = 1,...,n}$, the metric is represented by a field of $n\cx n$ positive symetric matrices $g_{ij}(x)$. The Laplacian can be expressed in these coordinates as an order 2 differential operator :
$$\Delta_g = -\frac{1}{\sqrt{g}}\sum_{i,j}\frac{\bd}{\bd x_i}g^{ij}\sqrt{g}\frac{\bd}{\bd x_j},$$
where $g^{ij}$ is the inverse of the matrix $g_{ij}(x)$ and $\sqrt{g}$ is the square-root of the determinant of $g_{ij}$ (see \cite{Cha84} p. 5 for this computation). 

A $\Cl C^1$ function $\phi$ satifies \ind{Dirichlet boundary conditions} if it vanishes on $\bd M$. If $\bd M$ is piecewise $\Cl C^1$, the function $\phi$ satisfies \ind{Neumann boundary conditions} if $\nabla \phi$ is tangent to $\bd M$ almost everywhere.
\pgh
We note $L^2(M)$ the set of real-valued square-summable functions on $M$, $\Cl H^1(M)$ the set of $L^2$-functions whose gradient in the distribution sense is a square-summable vector field, and $\tilde{\Cl H}^2\subset \Cl H^1$ the set of functions whose Laplacian (in the distribution sense) is $L^2$. The Laplacian $\Delta_g$ with Neumann (or Dirichlet) boundary conditions extends to a non-compact operator on $L^2$, whose domain is the maximal subset $\Cl D(\Delta)\subset \tilde{\Cl H}^2$ on which the corresponding boundary conditions are satisfied. When the boundary is empty, the maximal domain of this extension is $\tilde{\Cl H}^2$.

\begin{Defi}
A real $\lambda$ is in the \ind{spectrum} of $\Delta_g$ with Neumann (resp. Dirichlet) boundary conditions if the operator $$\Delta_g : \Cl D(\Delta_g)\ra L^2$$ is not invertible. It is in the \ind{essential spectrum} of $\Delta_g$ if there exists a sequence of functions $\psi_n\in\Cl D(\Delta_g)$, orthogonal in $L^2$, such that 
$$\lim_{n\ra\infty} \norm{\Delta_g\psi_n - \lambda\psi_n}_{L^2(M)} = 0.$$
\end{Defi}
We will note $Spec(g)$ and $Spec^e(g)$ the spectrum of $\Delta_g$ and its essential spectrum ; $\lambda_0(g)$ and $\lambda_0^{ess}(g)$ will denote their infima. We have $Spec^e(g)\subset Spec(g)\subset\Bb R_+$.

\begin{Prop}
The \ind{discrete spectrum} of $\Delta_g$, defined by $Spec^d(g)=Spec(g)\bs Spec^e(g)$, is a set of discrete eigenvalues with finite multiplicities.
\end{Prop}
See for example \cite{Eich07}, p11.

\begin{Ex}
It is shown in \cite{LaxPhi82} that the essential spectrum of a geometrically finite non-compact hyperbolic $n$-manifold is the half line $\left[\frac{(n-1)^2}{4},\infty\right)$. Therefore, the discrete spectrum is non-empty if and only if $\lambda_0<\frac{(n-1)^2}{4}$. This occurs, for instance, for hyperbolic surfaces whose convex core has a boundary pinched enough, or for geometrically finite acylindrical hyperbolic $3$-manifolds (see \cite{CanMinTay99}). Various other examples of non-compact manifolds with non-empty discrete spectrum are given in \cite{FisHam05}, Section 4.
\end{Ex}

We assume for the rest of this section that $M$ is complete without boundary and that $\Delta_g$ has at least one eigenvalue $\lambda_0<\lambda_0^{ess}$ (otherwise Theorem 1 is empty).
\pgh
We define now the Sobolev spaces on our non-compact manifold, whose use will be crucial for our arguments. Our presentation relies on the article \cite{Aub76} of T. Aubin.

Let $(M,g)$ be a non-compact riemannian manifold. We fix a $\Cl C^{\infty}$ atlas $\g U = (U_a,\phi_a)_{a\in A}$ of $M$, which is assumed to be \ind{locally finite} :  for all $a\in A$, $\bar{U_a}$ is compact and intersects only a finite number of  $U_b$, $b\neq a$. Let $k>0$ be fixed, assume the metric $g$ is $\Cl C^{k+2}$ (for this atlas).
\pgh
For any differential $p$-form $\alpha : (TM)^p\ra\Bb R$, with $p\leq k+2$, and for all chart-coordinates $(x_i)_{i=1..n}$ around $x\in M$, we set 
\beq\label{eq:DefNormForm}
|\alpha|_g^2(x) = \sum_{i_1,...,i_p,j_1,...,j_p = 1}^n g^{i_1j_1}...g^{i_pj_p}\alpha_{i_1...i_p}\alpha_{j_1...j_p},\eeq
where $g^{ij}$ is the inverse matrix of the metric $g_{ij}(x)$.
Then, by definition, for any $\phi\in\Cl C^1(M,\Bb R)$,
$$|d\phi|_g^2(x) = g_x(\nabla\phi, \nabla\phi) =: |\nabla\phi|^2_g.$$
Let $\phi\in\Cl C^q(M,\Bb R),$ with $2\leq q\leq k$. We note $\nabla^0\phi = \phi$, and for any $2\leq p\leq q$,
$$\nabla^p\phi = \nabla^{p-1}d\phi : (TM)^p\ra \Bb R$$
the $p-1$-th covariant derivative of $d\phi$. It is a linear $p$-form. 

For any $q\leq k+2$, the \ind{Sobolev norm of order $q$} of the function $\phi$ is then defined by :
$$\norm{\phi}^2_{\Cl H^q(g)} = \sum_{p = 0}^q \int_M |\nabla_{g}^p\phi|_{g}^2(x)dv_{g}(x)\in\Bb R_+\cup\{\infty\}.$$
We note
$$\Cl C^{\infty,q}(g) = \left\{\phi\in\Cl C^\infty(M,\Bb R) ; \norm{\phi}_{\Cl H^q(g)}<\infty\right\}.$$

\begin{Defi}\label{def:Sobolev}
The \ind{Sobolev space} of ordre $q$ is the completion of $\Cl C^{\infty,q}(g)$ pour la norme $\norm{.}_{\Cl H^q(g)}$. We note this space $\Cl H^q(g)$, it is a Hilbert space for the norm $\norm{.}_{\Cl H^q(g)}$.
\end{Defi}

Obviously, for any $p\geq q\geq 0$, $\Cl H^p\subset \Cl H^q$. Let us note that the expression of $\norm{.}_{\Cl H^q(g)}$ includes the partial derivatives of the metric $g$ up to order $q-1$. By definition, $H^0(g_0) = L^2(M,g_0)$, and $\Cl H^1(M)$ is canonically identified with the set of $L^2$ functions whose gradient in the distribution sense is an $L^2$ vector field, which we gave previously as a definition of $\Cl H^1$. Recall that we call $\tilde{\Cl H}^2(g)$ the set of $\Cl H^1$ functions whose Laplacian is a square summable function. 

\begin{Prop}
If $(M,g)$ is complete without boundary, with a lower bound for the Ricci tensor and positive injectivity radius, then 
$$\tilde{\Cl H}^2(g) = \Cl H^2(g),$$
and the corresponding norms are équivalent.
\end{Prop}
A proof of this result can be found in \cite{Heb96}, p16. Without any geometric assumption on $M$, for any $\Cl C^2$ function $u : M\ra \Bb R$, 
$$\left| \Delta_g u\right |^2 \leq n\left| \nabla^2 u\right |^2.$$
Therefore, $\Cl H^2(g)\subset\tilde{\Cl H}^2(g),$ and the Laplacian always maps $\Cl H^2$ into $L^2 = \Cl H^0$. We will need later to extend that property to Sobolev spaces of order $q>2$ ; this will require some more hypotheses on the curvature tensor.

\begin{Defi}
We say $(M,g)$ has \ind{bounded geometry at order $0$} if and only if its injectivity radius is positive and its sectional curvature is uniformly bounded on $M$. It has \ind{bounded geometry at order $k\geq 1$} if it has bounded geometry at order $0$ and there exists $C>0$ such that for all $p\leq k$, and for all $x\in M$,
$$\left|\nabla_{g}^p\Rm_{g}\right|(x)\leq C,$$
where $Rm_g$ is the Riemann curvature tensor of $g$.
\end{Defi}
As the Riemann tensor is expressed using derivatives of the metric up to order $2$, bounded geometry at order $k$ is only defined for $\Cl C^{k+2}$ metrics. For any $q\in\Bb N$, let us define a new norm on $\Cl C^\infty(M,\Bb R)$ by
$$\norm{\phi}^2_{\tilde{\Cl H}^q(g)} = \sum_{p = 0}^{q/2} \int_M |\Delta_{g}^p\phi|^2(x)dv_{g}(x)+\sum_{p = 0}^{q/2-1} \int_M |\nabla_g\Delta_{g}^p\phi|^2(x)dv_{g}(x)\in\Bb R_+\cup\{\infty\}$$
when $q$ is even, and
$$\norm{\phi}^2_{\tilde{\Cl H}^q(g)} = \sum_{p = 0}^{(q-1)/2} \int_M |\Delta_{g}^p\phi|^2(x)dv_{g}(x)+\sum_{p = 0}^{(q+1)/2} \int_M |\nabla_g\Delta_{g}^p\phi|^2(x)dv_{g}(x)\in\Bb R_+\cup\{\infty\}$$
when $q$ is odd. We note
$$\tilde{\Cl C}^{\infty,q}(g) = \left\{\phi\in\Cl C^\infty(M,\Bb R) ; \norm{\phi}_{\tilde{\Cl H}^q(g)}<\infty\right\},$$
and $\tilde{\Cl H}^q(g)$ the completion of $\tilde{\Cl C}^{\infty,q}(g)$ pour la norme $\norm{.}^2_{\tilde{\Cl H}^q(g)}$. By defintion,  $\tilde{\Cl H}^0(g) = L^2(g)$, $\tilde{\Cl H}^1(g) = \Cl H^1(g)$, and this definition of $\tilde{\Cl H}^2(g)$ is coherent with the one we gave in the previous paragraph.

\begin{Theo}\label{theo:LapSobolev}
Let $(M,g)$ be a riemannian complete manifold without boundary, and $k\in\Bb N$ such that $g$ is $\Cl C^{k+2}$ and $(M,g)$ has bounded geometry at order $k$. Then for any $q\leq k+2$, the set of $\Cl C^\infty$ functions with compact support is dense in $\Cl H^q(g)$, $\tilde{\Cl H}^q(g) = \Cl H^q(g),$ and the corresponding norms $\norm{.}_{\tilde{\Cl H}^q(g_0)}$ and $\norm{.}_{\Cl H^q(g_0)}$ are équivalent.
\end{Theo}

Our first assertion is given by the Theorem 2 of \cite{Aub76}. The second follows from Proposition 3 of the same paper. The following consequences comes immediately from the definition of the $\tilde{\Cl H}^q(g)$:

\begin{Coro}\label{coro:LapSobolev}
Under the same hypotheses, the Laplacian $\Delta_{g}$ is a bounded operator from $\Cl H^{k+2}(g)$ into $\Cl H^{k}(g_0)$, all its $L^2$-eigenfunctions are in $\Cl H^{k+2}(g)$ and for any $\lambda\in\Bb R$ which is not in the $L^2$-spectrum of $\Delta_{g}$, the operator $\Delta_{g}-\lambda$ is invertible from $\Cl H^{k+2}(g)$ onto $\Cl H^{k}(g)$.
\end{Coro}

\begin{Rema}
For these two last theorems, we may not need the assumption of a positive injectivity radius ; see for example the artcile by G. Salomonsen, \cite{Sal01}. Similarly, using harmonic coordinates as in Chapter 2 of \cite{Heb96} should allow to control only the Ricci tensor, and not the whole Riemann tensor. We will limit our presentation to manifolds with bounded geometry, for in that case the work of T. Aubin gives immediately the Theorem \ref{theo:LapSobolev}.
\end{Rema}

\subsection{Strong topology on complete manifolds}

We present now a topology on the set of metrics on $M$, the \ind{strong topology}, also known as the \ind{Whitney topology}. Even if it is classical in functional analysis on non-compact spaces (see for example \cite{Hir94}), it seems that this topology has never been used to study the properties of elliptic operators on non-compact manifolds.
\pgh
Recall all our manifolds are assumed to be $\sigma$-compact, and let $\g U = (U_a,\psi_a)_{a\in A}$ be a locally finite $\Cl C^\infty$ atlas on $M$ (which always exists). All smooth properties of functions and metrics on $M$ will be taken relatively to the smooth structure given by $\g U$. Let $k\in\Bb N$, we note $\Cl G^k(M)$ the set of $\Cl C^k$ riemannian metrics on $M$.

For any $\Cl C^\infty$ function $\epsilon : M\ra \Bb R^*_+$ and any metric $g\in\Cl G^k(M)$, we set
$$U_{k,\epsilon}(g) = \left\{h\in\Cl G^k(M) ; \max_{a\: :\: x\in U_a}\sup_{0\leq p \leq k}\sup_{(\alpha_1,...,\alpha_p)\in[1,n]}\left|\frac{\bd^p(h - g)}{\bd x_{\alpha_1}...\bd x_{\alpha_p}}\right|(x)\leq \epsilon(x)\right\},$$
where for all $a\in A$, the $x_{\alpha_i}$ are the local coordinates defined by the chart $(U_a,\psi_a)$. Hence, a metric $h\in U_{k,\epsilon}(g)$ if and only if in any point, in the local charts defined by $\g U$, all its partial derivatives of order less or equal to $k$ are $\epsilon$-closed to those of $g$. 

\begin{Defi}\label{def:TopoForte}
We call \ind{$\Cl C^k$-strong topology} on $\Cl G^k(M)$, or \ind{Whitney topology} of order $k$, the topology generated by the $U_{k,\epsilon}(g)$, where $\epsilon$ runs over all smooth positive functions on $M$ and $g$ over all $\Cl C^k$ riemannian metrics on $M$.
\end{Defi}

If $\g U'$ is another locally finite atlas defining the same smooth structure on $M$ as $\g U$, the strong topology defined using $\g U'$-charts is identical to the previous one. 
We can of course replace in this definition $\Cl G^k(M)$ by any space of $\Cl C^k$ functions from $X$ to $\Bb R$, where $X$ is a smooth manifold (here, $X = T^2M$). Recall there are two other classical topologies on such a function space : the \ind{uniform topology} of order $k$, and the \ind{compact-open topology} of order $k$. When $M$ is compact, these three topologies are the same, generated by the norm
$$\norm{f}_{\infty,k} = \sup_{a\in A} \sup_{x\in U_a}\sup_{0\leq p\leq k}\sup_{(\alpha_1,...,\alpha_p)\in[1,n]}\left|\frac{\bd^pf}{\bd x_{\alpha_1}...\bd x_{\alpha_p}}\right|(x).$$

For a description of general properties of these three topologies, the reader may consult \cite{Hir94}. Let us point out that when $M$ is non-compact, the $\Cl C^k$-strong topology is not metrizable, on the opposite to the uniform topology and the compact-open topology. For any compact $K\subset M$, we will call \emph{uniform norm of order $k$} the norm defined on $\Cl C^k$ functions on $K$ by
\beq\label{eq:NormInfty}
\norm{f}_{\infty,K,k} = \sup_{0\leq p\leq k}\norm{\nabla^pf}_{\infty,K}= \sup_{0\leq p\leq k}\sup_{x\in K} |\nabla^p f|(x).\eeq

\pgh
The $\Cl C^k$-strong topology is finer that the uniform topology. Here is a fundamental property of its converging sequences :

\begin{Prop}\label{prop:ConvForte}
Let $(g_n)\in\Cl C^k(M)^{\Bb N}$ be a sequence converging to $g_\infty$ for the $\Cl C^k$-strong topology. There exists $N>0$ and a compact $K\subset M$ such that for any $n\geq N$,
$$g_n|_{M\bs K} = g_\infty|_{M\bs K}.$$
Moreover, for any $p=0,...,k,$
$$\norm{g_n-g_\infty}_{\infty,K,p}\ra 0.$$
\end{Prop}
Note that the compact $K$ generally depends on the sequence $(g_n)$.

\begin{proof}
Let $(K_p)_{p\in\Bb N}$ be a sequence of compact sets such that $K_p\subset \inter{K}_{p+1}$ and 
$$\bigcup_{p\geq 0} K_p = M,$$ let us show there exists $p,N>0$ such that for any $n\geq N$, 
$$g_n|_{M\bs K_p} = g_\infty|_{M\bs K_p}.$$ By contradiction, assume that for any $N,p\in\Bb N$, there is $n_{N,p}\geq N$ and $x_{N,p}\in M\bs K_p$, such that $$|g_{n_{N,p}}(x_{N,p})-g_\infty(x_{N,p})|\geq \epsilon_{N,p}>0.$$
We can assume without loss of generality that there exists a $\delta>0$ such that for any $N,p>0$, all balls center in $x_{N,p}$ with radius $\delta$ are disjoint. Let $\epsilon : M\ra\Bb R^*_+$ be a function such that for any $N\in\Bb N$, $\epsilon(x_{N,N})<\epsilon_{N,N}.$
As all the $x_{N,N}$ are isolated, there exists such a function. We note 
$$U_{k,\epsilon}(g_\infty) = \left\{h\in\Cl G^k(M) ; \sup_{0\leq p \leq k}|\nabla_{g_\infty}^p(h - g_\infty)(x)|\leq \epsilon(x)\right\},$$
by definition of the $\Cl C^k$-strong topology, it is a neighbourhood of $g_\infty$. One can immediately check that for any $N\geq 0$, $g_{n_{N,N}}\notin U_{k,\epsilon},$ contradicting the hypothesis that $(g_n)$ converges to $g_\infty$ for the $\Cl C^k$-strong topology. Therefore, there exists a compact $K\subset M$ such that $g_n|_{M\bs K} = g_\infty|_{M\bs K}.$ By definition of the $\Cl C^k$-strong topology, we immediately get for any $p=0,...,k,$
$$\norm{g_n-g_\infty}_{\infty,K,p}\ra 0.$$
\end{proof}

\begin{Rema}
The converse is obviously true : if $(g_n)$ is a sequence of metrics such that there is a compact $K$ and an integer $N>0$, such that for any $n\geq N$, 
$$g_n|_{M\bs K} = g_\infty|_{M\bs K},$$
and if the $g_n|_{K}$ converge uniformly (with all their partial derivatives of order less or equal to $k$) to $g_\infty|_{K}$, then the sequence $(g_n)$ converges to $g_\infty$ for the $\Cl C^k$-strong topology.
\end{Rema}

Bounded topology is an open property for the strong topology :

\begin{Prop}\label{prop:StabBound}
Let $k\geq 2$ and $g_0\in\Cl G^k(M)$ be a metric with bounded geometry at order $k-2$. Then there exists a neighbourhood $\Cl U(g_0)$ for the $\Cl C^k$-strong topology such that any metric $g\in\Cl U(g_0)$ has bounded geometry at order $k-2$. Moreover, we can choose $\Cl U(g_0)$ such that for any $g\in \Cl U(g_0)$, the injectivity radius of $(M,g)$ has a uniform lower bound on $\Cl U(g_0)$, and the bound on the derivatives of $Rm_g$ is uniform on $\Cl U(g_0)$.
\end{Prop}

\begin{proof}[Sketch of proof]
The Riemann tensor is a polynomial into the partial derivatives of the metric of order up to $2$ with fixed coefficients. Therefore, a neighbourhood with uniform bounds on Riemann tensor's covariant derivatives of order $\leq k-2$ comes as a straightforward consequence of the definition of the $\Cl C^k$-strong topology. We let the reader check that if the injectivity radius $Inj(M,g_0)=\delta>0$, then for any $\epsilon>0$, there is a neighbourhood $\Cl U_\epsilon(g_0)$ of $g_0$ for the $\Cl C^k$-strong topology such that for any $g\in\Cl U_\epsilon(g_0), Inj(M,g)\geq \delta-\epsilon$.
\end{proof}

We will need the fact that two metrics close for the $\Cl C^k$-strong topology will define the same Sobolev spaces :

\begin{Prop}\label{prop:VoisFortSobolev}
For any $g\in\Cl G^k(M)$, any $p\in[0,k+1]$ and any constant $C>1$, there exists a neighbourhood $\Cl U_p$ of $g$ for the $\Cl C^k$-strong topology such that for any $h\in\Cl U_p$, and any function $f\in\Cl H^p(g)$,
$$\frac{1}{C}\int_M |\nabla_g^pf|_{g}^2dv_{g}\leq \int_M |\nabla_h^pf|_{h}^2dv_{h}\leq C\int_M |\nabla_g^pf|_{g}^2dv_{g},$$
these inequalities being always defined $\bar{\Bb R}$. In particular, for any $h\in\Cl U_p$, we have
$\Cl H^p(h) = \Cl H^p(g),$
and the norms $\norm{.}_{\Cl H^p(g)}$ et $\norm{.}_{\Cl H^p(h)}$ are uniformly equivalent.
\end{Prop}

\begin{proof}
Let $g\in\Cl G^k(M)$ and $C>1$. For any $h\in\Cl G^k(M)$ and any function $\phi : M\ra\Bb R+$, we have
$$\int_M \phi(x) dv_h(x) = \int_M \phi(x) \frac{dv_{h}}{dv_{g}}(x)dv_g(x).$$
For $p = 0$, let $1<C_0<C$, and let us consider the neighbourhood $\Cl U_0 = U_{k,\epsilon_0}(g)$ (see Definition \ref{def:TopoForte}) where $\epsilon_0 : M\ra\Bb R_+^*$ is such that for any $h\in U_{k,\epsilon_0}(g)$, and for any $x\in M$,
$$\frac{1}{C_0}\leq\frac{dv_{h}}{dv_{g}}(x)\leq C_0.$$ As $\frac{dv_{h}}{dv_{g}}(x)$ is the determinant of $h$ in a $g$-orthonormal basis of $T_xM$, such a function exists by continuity of the determinant. One can immediately check that $\Cl U_0$ satisfies the conclusion of our proposition at order $0$.

Let $p\in[1,k+1]$, we will use the following lemma :

\begin{Lemm}\label{lem:SoboCompact}
For any chart $U_a\in \g U$ with compact closure $K$, and for any (constant) $\eta>0$, there is a constant $C_{\eta,a}>1$ such that if $\norm{g-h}_{\infty,K,p}\leq \eta,$ then for any point $x\in U_\alpha$, and any function $f\in \Cl H^p(g)$,
$$\frac{1}{C_{\eta,a}}|\nabla_g^pf|_{g}^2(x)\leq |\nabla_h^pf|_{h}^2(x)\leq C_{\eta,a}|\nabla_g^pf|_{g}^2(x).$$
Moreover, we can take $C_{\eta,a}$ such that
$$\lim_{\eta\ra 0}C_{\eta,a} = 1.$$
\end{Lemm}
\begin{proof}
Let $a\in A$, set $K = \bar{U_a}$. We start with $p = 1$ : for any $x\in K$, let $g_{ij}(x)$ note the metric $g$ expressed in the coordinates defined by the chart $(U_a,\psi_a)$, and $g^{ij}(x) = \left(g_{ij}(x)\right)^{-1}$. By continuity of the matrix inverse, and uniform continuity of $g$ on $K$, for any $\eta>0$, there is a constant $C^1_{\eta,a}$ such that
$$\abs{g_{ij}-h_{ij}(x)}\leq \eta \Rightarrow \frac{1}{C^1_{\eta,a}}g^{ij}\leq h^{ij}\leq C^1_{\eta,a} g^{ij},$$
with moreover $$\lim_{\eta\ra 0}C^1_{\eta,a} = 1.$$

Let $p>1$ and $h\in \Cl G^k(M)$. On $U_a$, we write explicitely $|\nabla^p_hf|_h^2$ in the chart-coordinates defined by $(U_a,\psi_a)$. Translating the Christoffel symbols as sums of partial derivatives of the metric (see \cite{GHL04}, §2.B) at each step of the derivation, one can show by induction that this covariant derivative is a polynomial, whose variables are the partial derivatives of order $\leq p$ of $f$ and the partial derivatives of order $\leq p-1$ of the metric $h$, and whose coefficients do not depend on these. Moreover, this polynomial do not have any term which does not depend on the metric $h$. Therefore, the existence of the constants $(C_{\eta,a})_{\eta>0}$ satisfying the conclusion of our lemma comes directly from the uniform continuity of the partial derivatives of order $\leq p-1$ of the metric on $K$.
\end{proof}

Let now $\zeta = (\zeta_a)_{a\in A}$ be a partition of unity subordinated to $\g U$ : for any $a\in A, \zeta_a$ has compact support, and for any $x\in M$,
$$\sum_{a\in A}\zeta_a(x) = 1.$$ Let $f\in \Cl H^p(g)$, for any metric $h\in \Cl U_1$, we have
$$\norm{f}_{\Cl H^p(h)} = \sum_a\int_{U_a}\zeta_a(x)\left|\nabla_h^pf\right|(x)dv_h(x).$$
Let $(\epsilon_a)_{a\in A}$ be a family of positive constants such that for any $a\in A$ and $x\in U_a$, $\epsilon_a\leq \epsilon_0(x)$, where $\epsilon_0$ was defined above for $p = 0$. Let $C_a = C_{\epsilon_a,a}$ be the family of constants given by the previous lemma. Taking the $\epsilon_a$ small enough, we can assume for all $a\in A$, $1<C_a.C_0\leq C$. Let $\Cl U_p$ be a neighbourhood of $g$ such that for any $h\in \Cl U_p$ and any $U_a\in\g U$, we have
$$\norm{g-h}_{\infty,U_a,p}<\epsilon_a.$$ Such a neighbourhood exists by definition of the $\Cl C^k$-strong topology. We get then
$$\sum_{a\in A}\int_{U_a}\zeta_a(x)\left|\nabla_h^pf\right|(x)dv_h(x)\leq \sum_{a\in A}\int_{U_a}\zeta_a(x)C_a\left|\nabla_g^pf\right|(x)C_0dv_g(x),$$
which gives
$$\norm{f}_{\Cl H^p(h)} \leq C\norm{f}_{\Cl H^p(g)}.$$
We show similarly that $\norm{f}_{\Cl H^p(h)}\geq \frac{1}{C}\norm{f}_{\Cl H^p(g)}$, which concludes the proof of the Proposition \ref{prop:VoisFortSobolev}.
\end{proof}

To finish this part, we show that this topology is adapted to the study of the spectrum of the Laplacian. We say the Laplacian $\Delta_g$ has $q+1$ \ind{first eigenvalues} $\lambda_0(g)<\lambda_1(g)\leq ...\leq \lambda_q(g)$ if $Spec^e(\Delta_g)\subset (\lambda_q(g),\infty)$ and the first eigenvalues of the discrete spectrum (repeated according to their mutliplicity) are $\lambda_0(g)<\lambda_1(g)\leq ...\leq \lambda_q(g)$. The following theorem shows that the bottom of the essential spectrum and the lower part of the discrete spectrum vary continuously for the strong topology. 

\begin{Theo}\label{th:StabTrou}
Let $k\geq 1$. The application $g\fa \lambda_0^{ess}(g)$ is continuous from $\Cl G^k(M)$ to $\Bb R$ for the $\Cl C^k$-strong topology. Assume that for some $g_0\in\Cl G^k(M)$ and the Laplacian $\Delta_{g_0}$ has $q+1$ first eigenvalues isolated from the essential spectrum :
$$\lambda_0(g_0)<\lambda_1(g_0)\leq ...\leq \lambda_q(g_0)<\lambda_0^{ess}.$$
Then there is a neighbourhood $\Cl U(g_0)$ of $g_0$ for the $\Cl C^k$-strong metric such that for any $g\in\Cl U(g_0)$, the Laplacian $\Delta_g$ has $q+1$ first eigenvalues. Moreover, for any $0\leq j\leq q+1$, the map $g\fa \lambda_j(g)$ is continuous from $\Cl U(g_0)$ to $\Bb R$ for the $\Cl C^k$-strong topology.
\end{Theo}

\begin{proof}
Let $k\geq 1$, we separate the proof of this theorem into the following lemmas.

\begin{Lemm}
The application $g\fa \lambda_0^{ess}(g)$ is continuous from $\Cl G^k(M)$ to $\Bb R$ for the $\Cl C^k$-strong topology.
\end{Lemm}

\begin{proof}
Let $g\in \Cl G^k(M)$ and $\epsilon>0$. By Proposition 2.1 of \cite{Don81}, for any $h\in\Cl G^k(M)$ and any $\lambda>0$, the interval $(-\infty,\lambda]$ intersects the essential spectrum of $\Delta_h$ if and only if for any $\delta>0$, there exists an infinite-dimensional subspace $\Cl D_\delta\subset \Cl H^2(h)$ such that for any $f\in\Cl D_\delta$,
$$\int_M f\Delta_hf dv_h<(\lambda+\delta)\int_M f^2dv_h,$$
i.e. for any $f\in\Cl D_\delta$,
$$\frac{\int_M\abs{\nabla f}_h^2dv_h}{\int_M f^2dv_h}<\lambda+\delta.$$
Let us write $\lambda_0^{ess} = \lambda_0^{ess}(g)$ and for any $\delta>0$, we set $\Cl D_\delta\subset \Cl H^2(g)$ an infinite dimensional subspace of $\Cl H^2(g)$ such that $f\in\Cl D_\delta$,
$$\frac{\int_M\abs{\nabla f}_g^2dv_g}{\int_M f^2dv_g}<\lambda_0^{ess}+\delta.$$

Let $\eta>0$ be such that $(1+\eta)^2\lambda_0^{ess}<\lambda_0^{ess}+\epsilon$ and $(1+\eta)^2(\lambda_0^{ess}-\epsilon)<\lambda_0^{ess}-\epsilon/2$, and let $\Cl V_\eta$ be a neighbourhood of $g$ on which the Sobolev spaces are $(1+\eta)$-uniformly equivalent, as given by the Proposition \ref{prop:VoisFortSobolev}. For any $h\in\Cl V_\eta$ and any $f\in\Cl D_\delta$, we get then
$$\frac{\int_M\abs{\nabla f}_h^2dv_h}{\int_M f^2dv_h}<(1+\eta)^2(\lambda_0^{ess}+\delta)<\lambda_0^{ess}+\epsilon+(1+\eta)^2\delta$$
by definition of $\eta$. Therefore, $\lambda_0^{ess}(h)<\lambda_0^{ess}+\epsilon$ for any $h\in\Cl V_\eta$. 
Moreover, for any function $f\in\Cl H^2(g)$, we also have
$$\frac{\int_M\abs{\nabla f}^2_gdv_g}{\int_M\nabla f^2dv_g}\leq (1+\eta)^2\frac{\int_M\abs{\nabla f}^2_hdv_h}{\int_M\nabla f^2dv_h}.$$
Assume by contradiction there exists $h\in\Cl V_\eta$ such that $\lambda_0^{ess}(h)<\lambda_0^{ess}-\epsilon.$ Then, for any $\delta>0$ there exists $\Cl D'_\delta\subset \Cl H^2(g)$ infinite-dimensional such that for any $f\in\Cl D'_\delta$, we get
$$\frac{\int_M\abs{\nabla f}_h^2dv_h}{\int_M f^2dv_h}<\lambda_0^{ess}-\epsilon-\delta.$$
Therefore by definition of $\eta$, for any $\delta>0$ and any $f\in\Cl D'_\delta$,
$$\frac{\int_M\abs{\nabla f}_g^2dv_g}{\int_M f^2dv_g}\leq (1+\eta)^2(\lambda_0^{ess}-\epsilon-\delta)<\lambda_0^{ess}-\frac{\epsilon}{2}-(1+\eta)^2\delta.$$
The bottom of the essential spectrum of $\Delta_g$ is hence less than $\lambda_0^{ess}-\epsilon/2$, a contradiction.

We have eventually shown that for any $g\in\Cl G^k(M)$ and $\epsilon>0$, there exists a neighbourhood $\Cl V_\epsilon$ of $g$ for the $\Cl C^k$-strong topology such that for any $h\in\Cl V_\eta$, $$\lambda_0^{ess}(h)\in(\lambda_0^{ess}(g)-\epsilon,\lambda_0^{ess}(g)+\epsilon).$$
This concludes the proof of our lemma.
\end{proof}

Let $g_0\in\Cl G^k(M)$ be such that the Laplacian $\Delta_{g_0}$ has $q+1$ first eigenvalues isolated from the essential spectrum : there exists $\epsilon>0$ such that
$$\lambda_0(g_0)<\lambda_1(g_0)\leq ...\leq \lambda_q(g_0)\leq\lambda_0^{ess}(g_0)-2\epsilon.$$ Let $\Cl U(g_0)$ be a neighbourhood of $g_0$ for the $\Cl C^k$-strong topology such that for any $g\in\Cl U(g_0), \lambda_0^{ess}(g)>\lambda_0^{ess}(g_0)-\epsilon.$ Such a neighbourhood exists by the previous lemma. Let $\eta>0$ such that $(1+\eta)\lambda_q(g_0)\leq \lambda_0^{ess}(g_0)-\epsilon.$ Up to reducing $\Cl U(g_0)$ to a smaller neighbourhood, by Proposition \ref{prop:VoisFortSobolev} we can assume that for any $g\in\Cl U(g_0)$, and any $f\in\Cl H^1(g_0)$, 
$$\frac{\norm{\nabla_g f}_{L^2(g)}^2}{\norm{f}_{L^2(g)}^2}\leq (1+\eta)\frac{\norm{\nabla_{g_0} f}_{L^2(g_0)}^2}{\norm{f}_{L^2(g_0)}^2}.$$

\begin{Lemm}\label{lem:ContiLq}
For any $g\in\Cl U(g_0)$, the Laplacian $\Delta_g$ has $q+1$ first eigenvalues and for any $j\leq q+1$, the map $g\fa \lambda_j(g)$ is continuous from $\Cl U(g_0)$ to $\Bb R$ for the $\Cl C^k$-strong topology.
\end{Lemm}

To show the continuity of these eigenvalues, we will need the following version of the Min-Max principle for non-compact manifolds.

\begin{Prop}[Weak Min-Max principle on non compact manifolds]
Let $(M,g)$ be a Riemannian manifold such that the bottom of the essential spectrum $\lambda_0^{ess}(\Delta_g)>0$. For any $\lambda\in[0,\lambda_0^{ess}(g))$ and any $j\in\Bb N$, assume there exists a $j+1$-dimensional subspace $\Cl D_j\subset \Cl H^1(g)$ such that for any $f\in\Cl D_j$,
$$\frac{\norm{\nabla_g f}_{L^2(g)}^2}{\norm{f}_{L^2(g)}^2}\leq \lambda.$$
Then $\Delta_g$ has $j+1$ first eigenvalues smaller or equal to $\lambda$.
\end{Prop}
This is a very classical result, whose proof is the same as in the compact case for we assume to be below the bottom of the essential spectrum. We shall not prove it here, see for example the first chapter of \cite{Cha84}.

\begin{proof}[Proof of Lemma \ref{lem:ContiLq}]
Let us first show that for any $g\in\Cl U(g_0)$, the Laplacian $\Delta_g$ has $q+1$ eigenvalues. For any $0\leq j\leq q$, we note $\phi_i^0$ the eigenfunction of $\Delta_{g_0}$ associated to the eigenvalue $\lambda_j(g_0)$. Let $\Cl D_{q}(g_0)$ be the subspace of $\Cl H^1$ generated by the $\{\phi_j^0 ; 0\leq j\leq q\}$. One can easily check that for any $f\in\Cl D_{q}(g_0),$
$$\frac{\norm{\nabla_{g_0} f}_{L^2(g_0)}^2}{\norm{f}_{L^2(g_0)}^2}\leq \lambda_q(g_0).$$
By construction of $\Cl U(g_0)$, we have then for any $f\in\Cl D_{q}(g_0)$,
$$\frac{\norm{\nabla_{g} f}_{L^2(g)}^2}{\norm{f}_{L^2(g)}^2}\leq(1+\eta)\frac{\norm{\nabla_{g_0} f}_{L^2(g_0)}^2}{\norm{f}_{L^2(g_0)}^2}\leq \lambda_0^{ess}(g_0)-\epsilon<\lambda_0^{ess}(g).$$
Therefore, by the previous Min-Max principle, $\Delta_g$ has $q+1$ first eigenvalues in $[0,\lambda_0^{ess}(g_0)-\epsilon]$.
\pgh
For any $0\leq j\leq q$, let us show $g\fa \lambda_j(g)$ is continuous for the $\Cl C^k$-strong topology. Let $g\in\Cl U(g_0)$, we note $\lambda_0<\lambda_1\leq...\leq \lambda_j$ the $j+1$ first eigenvalues of $g$, $\phi_0,...,\phi_j$ the associated eigenfunctions, and $\Cl D_j(g)$ the $j+1$-dimensional space they generate. Let fix $\epsilon>0$, and let $\eta>0$ be such that $$(1+\eta)\lambda_j\leq \lambda_j+\epsilon \mbox{ and } (1+\eta)(\lambda_j-\epsilon)\leq \lambda_j-\frac{\epsilon}{2}.$$
Let $\Cl U_\epsilon(g)\subset \Cl U(g_0)$ be a neighbourhood of $g$ such that for any $h\in\Cl U_\epsilon(g)$ and any function of $\Cl H^1(g)$,
$$\frac{1}{1+\eta}\frac{\norm{\nabla_{g} f}_{L^2(g)}^2}{\norm{f}_{L^2(g)}^2}\leq \frac{\norm{\nabla_{h} f}_{L^2(h)}^2}{\norm{f}_{L^2(h)}^2}\leq(1+\eta)\frac{\norm{\nabla_{g} f}_{L^2(g)}^2}{\norm{f}_{L^2(g)}^2}.$$
Such a neighbourhood exists by Proposition \ref{prop:VoisFortSobolev}. We claim that for any $h\in \Cl U_\epsilon(g)$, $\lambda_j(h)\in[\lambda_j(g)-\epsilon, \lambda_j(g)+\epsilon]$. Indeed, for any $f\in\Cl D_j(g)$, we have
$$\frac{\norm{\nabla_{h} f}_{L^2(h)}^2}{\norm{f}_{L^2(h)}^2}\leq(1+\eta)\frac{\norm{\nabla_{g} f}_{L^2(g)}^2}{\norm{f}_{L^2(g)}^2}\leq (1+\eta)\lambda_j\leq \lambda_j+\epsilon,$$
which implies by the Min-Max principle that $\lambda_j(h)\leq \lambda_j+\epsilon.$ Assume now by contradiction there is $h\in\Cl U_\epsilon(g)$ such that $\lambda_j(h)<\lambda_j-\epsilon$. Let $\Cl D_j(h)$ be the space generated by the $j+1$ first eigenfunctions of $\Delta_h$. For any $f\in\Cl D_j(h)$, we have
$$\frac{\norm{\nabla_{g} f}_{L^2(g)}^2}{\norm{f}_{L^2(g)}^2}\leq(1+\eta)\frac{\norm{\nabla_{h} f}_{L^2(h)}^2}{\norm{f}_{L^2(h)}^2}\leq (1+\eta)(\lambda_j-\epsilon)\leq \lambda_j-\frac{\epsilon}{2}.$$
This implies that $\lambda_j(g)\leq \lambda_j-\frac{\epsilon}{2}<\lambda_j(g)$, a contradiction.
\end{proof}

\begin{Rema}
The continuity of the first eigenvalues for the strong topology could also be deduced from the continuity of the spectrum for Lipschitz convergence, which is proved in Theorem A. of \cite{FisHam05}.
\end{Rema}

\subsection{Proof of Genericity Theorem}
Let $(M^n, g_0)$ be a riemannian manifold with bounded geometry at order $k\geq\frac{n}{2}$ such that the Laplacian $\Delta_{g_0}$ has $q+1$ first eigenvalues $\lambda_0<\lambda_1\leq...\leq\lambda_q$ (repeated according to their multiplicity) isolated from the essential spectrum of $M$. By the results of our previous paragraph, we know there is a neighbourhood $\Cl U(g_0)$ for the $\Cl C^{k+2}$-strong topology such that for any metric in $\Cl U(g_0)$, the geometry is bounded at order $k$, and the Laplacian has $q+1$ first eigenvalues isolated from the essential spectrum. To finish the proof of our Theorem \ref{th:PhiM}, we have now to show there is a generic set $\Cl U_M(g_0)\subset \Cl U(g_0)$ such that for any metric in $g\in\Cl U_M(g_0)$, the $q+1$ first eigenvalues of $\Delta_g$ are simple and the associated eigenfunctions are Morse.

\begin{Defi}
Let $E$ be a topological space. A subset $A\subset E$ is a \ind{generic set} in $A$ if it can be obtained as a countable intersection of dense open sets in $E$. 
\end{Defi}
Such a subset has also been called a \ind{Baire set}, a \ind{residual set}, a \ind{dense $G_\delta$} or a \ind{set of second category}. 

Baire's theorem claims that a generic subset of a \emph{complete metric space} E is dense in $E$. The $\Cl C^k$-strong topology we defined in the previous paragraph is not metrizable. Even though Baire's theorem do not apply, a similar statement is still true :

\begin{Theo}\label{th:GBaire}
The set $\Cl G^{k+2}(M)$ with the $\Cl C^{k+2}$-strong topology is a \ind{Baire space}, which means that any generic subset of $\Cl G^{k+2}(M)$ is dense.
\end{Theo}

This is a direct corollary of Theorem 4.2 p 59 of \cite{Hir94} and of the fact that if a sequence of metrics $(g_n)$ uniformly converges on $M$ with its partial derivatives up to order $k$, then the limit is still a $\Cl C^k$ metric on $M$.

\paragraph{Simple eigenvalues are generic}

We note $\Cl G^{k+2}_q(M)$ the set of $\Cl C^{k+2}$-metrics on $M$ $g$ with bounded geometry at order $k$, whose Laplacien admits $q+1$ first eigenvalues $\lambda_0(g)<\lambda_1(g)\leq ...\leq \lambda_q(g)<\lambda_0^{ess}(g).$ By Proposition \ref{prop:StabBound} and Theorem \ref{th:StabTrou}, $\Cl G^{k+2}_q(M)$ is open in $\Cl G^{k+2}(M)$. Our first genericity result is the following :

\begin{Theo}\label{th:GeneSimpleEigval}
Given any metric $g\in\Cl G^{k+2}_q(M)$ and any non-empty open set $O\subset M$, there exists a metric $h\in\Cl G^{k+2}_q(M)$, equal to $g$ outside $O$, and arbitrarily closed to $g$ for the $\Cl C^{k+2}$-strong topology on $M$, such that all eigenvalues of $\Delta_h$ are distinct.
\end{Theo}

\begin{Coro}
The set of metrics $g\in\Cl G^{k+2}_q(M)$ such that the $q+1$ first eigenvalues of $\Delta_g$ are distinct is open and dense in $\Cl G^{k+2}_q(M)$.
\end{Coro}

\begin{proof}[Proof of the corollary]
The openness of the set of metrics with distinct $q+1$ first eigenvalues comes from Theorem \ref{th:StabTrou}. Its density comes from the previous theorem.
\end{proof}

Our proof of this theorem is inspired from the proof of Theorem 1 of K. Uhlenbeck's paper \cite{Uhl76} which we adapt to the non-compact case.

\begin{proof}[Proof of Theorem \ref{th:GeneSimpleEigval}]
The key tools for the proof of this genericity statement, as well as other generic properties we will prove, are Fredholm maps and Sard-Smale Transversality Theorem we present now. 

\begin{Defi}
\bi
\item An application $f : E\ra G$ between two Banach-spaces is a \ind{Fredholm operator} if its kernel is finite dimensional, its image is closed, and its co-kernel $G/Im(F)$ is finite dimensional. 
\item The \ind{index} of a Fredholm operator $F$ is $dim(Ker F)-dim(CoKer F)$. 
\item A differentiable application $F : \Cl X\ra \Cl Y$ between two Banach manifolds is a \ind{Fredholm map} if its differential is a Fredholm operator in each point. Its index is then, in each point, the index of its differential.
\ei
\end{Defi}

The concept of Fredholm map was introduced by S. Smale in \cite{Sma65}. He shows that when the basis manifold is connected, the index does not depend on the chosen point.

\begin{Ex}
Let $g\in\Cl G^{k+2}_q(M)$. For any $\lambda\in (0,\lambda_0^{ess}(g_0))$, the operator $\Delta_{g}-\lambda$ is a Fredholm operator from $\Cl H^{k+2}$ to $\Cl H^{k}$, for $\Delta_{g}$ is elliptic and its eigenvalues lower than $\lambda_0^{ess}(g_0)$ have finite multiplicities (see \cite{Sma65}, Section 2, and references given there). It has index 0 because it is self-adjoint.
\end{Ex}

Recall a \emph{regular value} of $F : \Cl X\ra \Cl Y$ is a point $y\in \Cl Y$ such that for any $x\in F^{-1}(y)$, $D_xF$ is surjective.

\begin{Theo}[Sard-Smale Theorem]\label{th:SardSmale}
Let $F : \Cl X\ra \Cl Y$ be a Fredholm map between separable Banach manifolds. If $F$ is $\Cl C^r$ with $r> index(F)$, then the set of regular values of $F$ is generic in $\Cl Y$.
\end{Theo}

See Smale's paper \cite{Sma65} for a proof. The following corollary of Sard-Smale Theorem is proved in \cite{Uhl76}, p1061 :

\begin{Coro}\label{coro:SardSmale}
Let $F : \Cl H\cx B\ra \Cl Y$ be a $\Cl C^r$ map, $\Cl H,\Cl B$ and $\Cl Y$ are Banach manifolds with $\Cl H$ and $\Cl Y$ separable. If for all $b\in \Cl B, F_b = F(.,b)$ is a Fredholm map of index $<r$, and if $0\in \Cl Y$ is a regular value for $F$ then the set 
$$\left\{b\in \Cl B, 0 \mbox{ is a regular value for }F_b\right\}$$
is generic in $\Cl B$.
\end{Coro}

Let $g\in\Cl G^{k+2}_q(M)$, we note $\lambda_0(g)<\lambda_1(g)\leq ...\leq \lambda_q(g)<\lambda_0^{ess}$ its first eigenvalues and $\phi_0^g,...,\phi_q^g$ the associated eigenfunctions. Let $\Cl U_g\subset \Cl G^{k+2}_q(M)$ be a neighbourhood of $g$ for the $\Cl C^{k+2}$-strong topology such that for any $h\in \Cl U_g$, $\lambda_g(h)<\lambda_0^{ess}(g)$ : such a neighbourhood exists by theorem \ref{th:StabTrou}. Let $O\subset M$ be a non-empty open set ; we can assume without generality that it has compact closure and smooth boundary. Let $K\subset O$ be a compact subset with non-empty interior and smooth boundary.

We note $\Cl G^{k+2}(O)$ the set of $\Cl C^{k+2}$ metrics on $\bar{O}$, and we define an injection$\iota_g : \Cl G^{k+2}(O)\ra \Cl G^{k+2}(M)$ the following way. Let $\chi_K : M\ra \Bb R$ be a smooth <<cut-off>> function such that $\chi_K|_K\equiv 1$ and $\chi_K|_{M\bs O} \equiv 0$. For all $h\in \Cl G^{k+2}(O)$, we define
$$\iota_g : h\fa \chi_K h + (1-\chi_K)g.$$
We note $\Cl O_g\subset \Cl G^{k+2}(O)$ the set of all $\Cl C^{k+2}$-metrics $h$ such that $\iota_g(h)\in \Cl U_g$, where $\Cl U_g$ is the neighbourhood of $g$ we defined previously. The set $\Cl O_g$ is open in $\Cl G^{k+2}(O)$ for the $\Cl C^{k+2}$-uniform topology on $\bar{O}$ : it is therefore a Banach manifold (getting a Banach manifold was the main purpose of restricting to compact perturbations of the metric). Remark that by the Decomposition Principle of \cite{DonLi79} we already mentionned, for all $h\in\Cl O_g$, 
$$\lambda_0^{ess}(\iota_g(h)) = \lambda_0^{ess}(g).$$

We can assume moreover, reducing $\Cl U_g$ if necessary, that for any $h\in\Cl O_g$ and any integer $g\leq k+2$, the Sobolev $\Cl H^q(M,g)$ and $\Cl H^q(M,\iota_g(h))$ are equivalent. We note from now on for all $q\in \Bb N, H^q = H^q(M,g)$: we still work avec the Sobolev space on the complete manifold.

We set $$\Cl S^{k+2} = \left\{u\in \Cl H^{k+2} ; \int_M{ u^2} dv_g= 1\right\} :$$ it is a Banach submanifold of $\Cl H^{k+2}$ with codimension $1$, and for any $\phi\in S^{k+2}$, the tangent space in $\phi$ to $\Cl S^{k+2}$ is
$$T_\phi \Cl S^{k+2} = \left\{v\in \Cl H^{k+2}(M), \int_M\phi v dv_g = 0\right\}.$$
From now on, for any $h\in\Cl O_g$, we will note $\Delta_h$ for $\Delta_{\iota_g(h)}$ and $\lambda_0(h)<\lambda_1(h)\leq...\leq \lambda_q(h)<\lambda_0(g)$ the $q+1$ first eigenvalues of $\Delta_h$.

Let us set now $$F : \Cl S^{k+2}\cx (0,\lambda_0^{ess}(g)) \cx \Cl O_g\ra \Cl H^{k}$$
defined $$F(\phi, \lambda, h) = \Delta_h\phi - \lambda\phi.$$ 

By the explicit formula of the Laplacian (see Paragraph \ref{ssec:BasicSpec}) and general properties of the evaluation map on a compact set (see \cite{AbRob67}, Chapitre 3), since for any  $\Cl C^2$-function $\phi : M\ra \Bb R$ the application 
$$(h,x)\in \Cl O_g\cx M\fa \Delta_h(\phi)(x)$$ is $\Cl C^{k}$, we know the application $F$ is $\Cl C^{k}$. We will note $$F_h = F(.,.,h).$$

\begin{Lemm}\label{lem:CaracValSimp}
The function $\phi\in S^{k+2}$ is an eigenfunction of $\Delta_h$ associated to $\lambda$ if and only if $F(\phi,\lambda,h) = 0$. Moreover, the eigenvalue $\lambda$ is simple if and only if $(\phi,\lambda)$ is a regular point of $F_h$.
\end{Lemm}
This lemma was already proved in \cite{Uhl76}, p1064.
\begin{proof}
The first assertion is obvious. Let $(\phi,\lambda,h)$ such that $\Delta_h\phi = \lambda\phi$, the differential $dF_h : T_\phi\Cl S^{k+2}\cx ]0,l[\ra \Delta_h(\Cl H^{k+2})$ is
$$dF_h(v,s) = \Delta_hv - s\phi - \lambda v,$$ where $v$ runs over all $\Cl H^{k+2}$ which are $L^2$-orthogonal to $\phi$ and $s\in \Bb R$ Moreover, $v\fa \Delta_hv - \lambda v$ is surjective on the $L^2$-orthogonal of $\phi$ in $\Cl H^k$ if and only if $\lambda$ is a simple eigenvalue. Varying parameter $s$, this concludes our proof.
\end{proof}

Our genericity result will be a direct consequence of the Corollary \ref{coro:SardSmale} of Sard-Smale Theorem and the following lemma :

\begin{Lemm}\label{lem:RegVal}
In any $(\phi,\lambda,h)\in F^{-1}(0)$, the differential $d_{(\phi,\lambda,h)}F$ is surjective : $0$ is a regular value of $F$.
\end{Lemm}

\begin{proof}
Let $(\phi,\lambda,h)\in F^{-1}(0)$. If $\lambda$ is a simple eigenvalue of $\Delta_h$ (which is the case when $\lambda = \lambda_0(h)$ by Courant's Nodal Theorem), then our lemma is a consequence of Lemma \ref{lem:CaracValSimp}. Therefore, we can assume $\lambda$ has multiplicity $m\in[2,q]$ : there are $\phi_2,...,\phi_m$ which are $L^2$-orthogonal to $\phi$ and to each other, such that for any $i\in [2,m], (\phi_i,\lambda,h)\in F^{-1}(0).$ We note $\bd_h F$ the partial derivative of $F$ along $T\Cl O_g$ : with our previous notations, the total derivatives splits into
$$d_{(\phi,\lambda,h)}F(v,s,\ddh) = dF_h(v,s) + \bd_hF(\ddh).$$
We want to show
$$Im(dF_h) \oplus Im(\bd_h F) = \Cl H^k.$$
We saw that $dF_h(v,s) = \Delta_hv - s\phi - \lambda v,$ therefore the $L^2$-orthogonal complement of $Im(dF_h)$ in $\Cl H^k$ is $Vect(\phi_2,...,\phi_m)$. To find an $m-1$-dimensional subspace $H'\subset Im(\bd_h F)$ such that $Im(dF_h) \oplus H' = \Cl H^k,$ it is enough to show that $Im(dF_h)+Im(\bd_h F)$ is $L^2$-dense in $\Cl H^k$. This will come from the following lemma :

\begin{Lemm}[Density of $Im(\bd_h F)$]\label{lemm:JDense}
Let $w\in\Cl L^2(M)$, $w\in\Cl C^2(M\bs\{y\})$ where $y\in M$. Assume that for any $j\in Im(\bd_h F)$,
$$\int_Mw(x)j(x)dv_h(x) = 0,$$
then $w$ is constant on $K$. 
\end{Lemm}

\begin{proof}
This is inspired from \cite{Uhl76}, p 1075. 
The tangent space to $\Cl O_g$ in the metric $h$ is
$$T_h\Cl O_g = \left\{\ddh ; \ddh \hbox{ symetric 2-form on } M \mbox{ with support in } O\right\}.$$

Let note $J = Im(\bd_h F)$, the fonction $w$ is $L^2$-orthogonal to $J$ if and only if for all $\ddh\in T_h\Cl O_g$, 
\beq\label{eq:OrthoJ}
\int_M w\bd_hF(\ddh)dv_h = 0.\eeq
As $$F(\phi,\lambda,h) = \Delta_h\phi-\lambda\phi,$$
we have by definition at a point $(\phi,\lambda,h)\in Q$, $$\bd_h F(\ddh) = \bd_h(\Delta_h\phi)(\ddh).$$
Moreover, by Green Formula we have
\beq\label{eq:GreenOrthog0}\int_M \iota_g(h)(\nabla_{\iota_g(h)}\phi,\nabla_{\iota_g(h)}w)dv_h = \int_M w(\Delta_h\phi)dv_h.\eeq
Let us note $d\phi\odot dw$ the symetric $2$-form defined for any $x\in M$ and $X,Y\in T_x M$ by
$$d\phi\odot dw(X,Y) = \frac{1}{2}(d_x\phi(X)d_xw(Y) + d_xw(X)d_x\phi(Y)).$$
In a any local coordinates basis $(x_i)$ on a neighbourhood of $x$, the matrix of $d\phi\odot dw$ is $\frac{1}{2}(d_{x_i}\phi d_{x_j}w + d_{x_j}\phi d_{x_i}w)$. As $\iota_g(h)$ is symetric, we have then at point $x$ and in this coordinates basis,
$$\iota_g(h)(\nabla_{\iota_g(h)}\phi,\nabla_{\iota_g(h)}w) = \sum_{i,j}(\iota_g(h))^{ij}d_{x_i}\phi\  d_{x_j}w = \Trace\left(\iota_g(h)(x)^{-1}(d_x\phi\odot d_xw)\right).$$
The Formula (\ref{eq:GreenOrthog0}) can then be rewritten as follows:
\beq\label{eq:GreenOrthog}
\int_M \Trace\left(\iota_g(h)^{-1}\cdot d\phi\odot dw\right)dv_h = \int_M w(\Delta_h\phi)dv_h,\eeq
where the inverse and the matrix product are taken in any coordinate basis in each point (the expression does not depend on that basis). 

Let us note
$$(T_h\Cl O_g)_0 = \left\{\ddh\in T_h\Cl O_g ; \forall x\in K', \Trace(h^T\ddh(x)) = 0\right\} :$$
$(T_h\Cl O_g)_0$ is the space of variations de $h$ which preserve the volume form $dv_h$ on $O$. It will be enough to consider such variations. One can immediately check that $d\iota_g(T_h\Cl O_g)_0$ is contained in the space of variations which preserve the volume form of $\iota_g(h)$ on all $M$, for the cut-off function only adds a multiplicative term to the volume element on $O\bs K$ which stays constant during these variations. 
Differentiating term by term the equation (\ref{eq:GreenOrthog}) along $(T_h\Cl O_g)_0$, we get then for any $\ddh\in (T_h\Cl O_g)_0$,
\beq\label{eq:FormOrthog}
\int_M\Trace\left(d(\iota_g)^{-1}(\ddh)\cdot dw\odot d\phi\right)dv_h = \int_M w\bd_hF(\ddh)dv_h = 0\eeq
by equation (\ref{eq:OrthoJ}). Moreover, for any $x\in K$, by definition we have
$$d\iota_g(\ddh)(x) = \ddh(x),$$
therefore
$$d(\iota_g)^{-1}(\ddh) = -g^{-1}\cdot\ddh.$$
By Formula (\ref{eq:FormOrthog}), we get then for allmost all $x\in K$ and all symetric matrix $\ddh$ such that $\Trace(h^{-1}\cdot\ddh) = 0$,
$$\Trace\left(-g^{-1}\cdot\ddh\cdot dw\odot d\phi\right) = 0.$$
Since
$$\langle A,B\rangle := \Trace(A\cdot B)$$ defines a scalar product on the set of symetric matrices, and as the application $\ddh\fa -g^{-1}\ddh$ is invertible, for almost all $x\in K$ there exists then $\lambda(x)\in\Bb R$ such that
$$dw\odot d\phi = \lambda(x)h.$$

Nevertheless, $h$ has maximal rank $n\geq 2$ and is positive. On the opposite, $dw\odot d\phi$ has rank at most 2, and if $n=2$, its determinant is negative. Therefore $$\lambda(x) = 0$$ almost everywhere $K$. As $dw$ is continuous outside $y$, its support is therefore included in the set of critical points of $\phi$ which are in $K$. As $M$ has infinite volume, $\phi$ is not constant. By the Unique Continuation principle of Aronszajn for solutions of elliptic equations (see \cite{Aro57}), the set of critical points of $\phi$ has empty interior. The differential $dw$ vanishes then on $K$, and  $w$ is constant.
\end{proof}

As $K$ has non-empty interior, by Aronszajn's Unique Continuation principle the $\phi_2,...,\phi_m$ cannot be constant $K$. Therefore $Im(dF_h)+Im(\bd_h F)$ is $L^2$-dense in $\Cl H^k$, which ends the proof of Lemma \ref{lem:RegVal}.

\end{proof}
Lemma \ref{lem:RegVal} and Corollary \ref{coro:SardSmale} implies that the set of metrics $h\in \Cl O_g$ such that the Laplacian of $\iota_g(h)$ has distinct eigenvalues is generic in $\Cl O_g$, hence dense, for the $\Cl C^{k+2}$-uniform topology. As $\iota_g$ is continuous from $\Cl G^{k+2}(O)$ with the $\Cl C^{k+2}$-uniform topology to $\Cl G^{k+2}(M)$ with the $\Cl C^{k+2}$-strong topology, this implies the theorem \ref{th:GeneSimpleEigval}.

\end{proof}

\paragraph{Morse eigenfunction are generic}

Recall a Morse function is a $\Cl C^2$-function whose Hessian is non-degenerate in all its critical points. For a $\Cl C^2$ function $f\in\Cl C^2(M,\Bb R)$, being Morse is well defined as long as a $\Cl C^2$-atlas is fixed on $M$. In particular, it does not depend on the metric. Morse functions are generic in $\Cl C^2(M)$ for the $\Cl C^2$-strong topology : this can be shown by an argument similar (but shorter) to the proof we give now that generic metrics have Morse eigenfunctions.

Let $g\in\Cl G^{k+2}_q(M)$, by Theorem \ref{th:GeneSimpleEigval} we can assume that the $q+1$ first eigenvalues of $\Delta_g$ are distinct. Let $\Cl U\subset\Cl G^{k+2}_q(M)$ be a neighbourhood of $g$ such that for any $h\in\Cl U$, the $q+1$ eigenvalues of $\Delta_h$ are distinct and the Sobolev Spaces $\Cl H^k(M,g)$ and $\Cl H^k(M,h)$ are equivalent. We will note $\Cl H^k = \Cl H^k(M,g)$, with its norm $\norm{.}_{\Cl H^{k+2}(g)}$. From now on, a \emph{<<compact $K\subset M$>>} will mean a compact subset of $M$ with smooth boundary and non-empty interior. Let us fix $i\in[0,q]$, for any compact $K\subset M$, we note $\Cl M_{\phi_i}(K)\subset \Cl U$ the set of metrics $h\in\Cl U$ such that the $i$-th eigenfunction $\phi_i^h$ is Morse on $K$. The first step of our main Theorem \ref{th:PhiM} is given by the following proposition:

\begin{Prop}\label{prop:KMorseOuvert}
For any compact $K\subset M$, $\Cl M_{\phi_i}(K)$ is an open subset of $\Cl U$ for the $\Cl C^{k+2}$-strong topology.
\end{Prop}

\begin{proof}
Let $h\in\Cl M_{\phi_i}(K)$, by hypothesis $\phi_i^h$ is Morse on $K$ and $h$ has bounded geometry at order $k$. By Corollary \ref{coro:LapSobolev}, the operator $\Delta_h$ is then bounded from $\Cl H^{k+2}$ to $\Cl H^{k}$, associated to the quadratic form
$$Q_h(f,\phi) = \int_M f(x)\Delta_h\phi(x)\alpha_h(x)dv_g(x),$$ where we have noted
$$\alpha_h(x) = \frac{dv_{h}}{dv_{g}}(x).$$
From the general expression of the Laplacian from the metric (see Paragraph \ref{ssec:BasicSpec}), the application $h\ra \Delta_h$ is continuous from $\Cl U$ with the $\Cl C^{k+2}$-strong topology to the set of bounded operators from $\Cl H^{k+2}$ to $\Cl H^{k}$. Therefore, the application which maps $h$ to the $i$-th eigenfunction $\phi_i^h$ of $\Delta_h$ is continuous from $\Cl U$ to $\Cl H^{k+2}$ (the reader used to spectral theory will find a proof of this result in the classical book T. Kato, \cite{Kat95}, Section IV.3.5 ; it can also be shown by a Implicit Function Theorem in $\Cl H^{k+2}$).

Now the set of Morse functions on $K$ is open for the $\Cl C^2$-uniform topology on $K$ : there exists $\eta>0$ such that any fonction $f : K\ra \Bb R$ which satisfies
$$\norm{f-\phi_i^g}_{\infty, K,2} := \sup_{0\leq p\leq 2}\sup_{x\in K}\abs{\nabla_g^p(f-\phi_i^g)}(x)<\eta$$ is Morse on $K$. 

Moreover, by Sobolev injections on the compact $K$ (see \cite{Heb96} p21), for $k\geq \frac{n}{2}$ (where $n$ is the dimension of $M$), $\Cl H^{k+2}(K)$ with the Sobolev norm restricted to $K$ injects continuously into $\Cl C^2(K)$ with the $\norm{.}_{\infty, K,2}$. There exists then $\delta>0$ such that for any fonction $f\in\Cl H^{k+2}(K)$ satisfying
$$\sup_{m\in[0,k+2]} \int_{K}\left|\nabla_g^m (f-\phi_i^g)\right|_{g}^2dv_{g}<\delta,$$
we get
$$\norm{f-\phi_i^g}_{\infty, K,2} <\eta$$

By continuity of the map $h\ra \phi_0^h$, there is a neighbourhood $\Cl U_1\subset \Cl U$ of $g$ for the $\Cl C^{k+2}$-strong topology such that for any $h\in\Cl U_1$, 
$$\sup_{m\in[0,k+2]} \int_{M}|\nabla_g^m (\phi_i^h-\phi_i^g)|^2dv_{g}<\delta.$$
By previous arguments, $\phi_i^h$ is Morse on $K$ for any $h\in\Cl U_1$, which concludes the proof of Proposition \ref{prop:KMorseOuvert}.
\end{proof}

The next crucial step of the proof of Theorem \ref{th:PhiM} will be the following proposition :

\begin{Prop}\label{prop:KMorseDense}
For any compact $K\subset M$, $\Cl M_{\phi_i}(K)$ is dense in $\Cl U$ for the $\Cl C^{k+2}$-strong topology.
\end{Prop}

Like for Theorem \ref{th:GeneSimpleEigval}, the proof of this proposition is inspired from \cite{Uhl76}, and relies on the use of an appropriate version of Sard-Smale Theorem.

\begin{proof}

Let $O\subset M$ be an open set with compact closure such that $K\subset O$, and $\chi_K : M\ra \Bb R$ be a cut-off function which is 1 on $K$ and $0$ outside $O$. We go back to the notations of the proof of Theorem \ref{th:GeneSimpleEigval}. $\Cl G^{k+2}(O)$ is the space of $\Cl C^{k+2}$ metrics on $\bar{O}$, $\iota_g : \Cl G^{k+2}(O)\ra \Cl G^{k+2}(M)$ is defined by 
$$\iota_g(h) = \chi_K h + (1-\chi_K)g.$$
The map $\iota_g$ is continuous from $\Cl G^{k+2}(O)$with the $\Cl C^{k+2}$-uniform topology to $\Cl G^{k+2}(M)$ with the $\Cl C^{k+2}$-strong topology.

Let $a,b\in\Bb R$ such that $\lambda_{i-1}(g)<a<\lambda_i(g)<b<\lambda_{i+1}(g),$ such $a,b$ exist by hypothesis. By continuity of $g'\fa \lambda_i(g')$ in $\Cl G^{k+2}(M)$ for the strong topology (see Theorem \ref{th:StabTrou}), there exists a neighbourhood $\Cl U_i\subset \Cl U$ such that for any $g'\in\Cl U_i$, $\lambda_{i-1}(g')<a<\lambda_i(g')<b<\lambda_{i+1}(g')$. We note $\Cl O_i\subset \Cl G^{k+2}(O)$ the set of all $\Cl C^{k+2}$-metrics $h$ such that $\iota_g(h)\in \Cl U_i$. For all $q\in \Bb N,$ we note $ H^q = H^q(M,g)$,
$$\Cl S^{k+2} = \left\{u\in \Cl H^{k+2} ; \int_M{ u^2} dv_g= 1\right\}.$$  and for any $h\in\Cl O_g$, we will note $\Delta_h$ for $\Delta_{\iota_g(h)}$, $\lambda_i(h)$ the $i$-th eigenvalue of $\Delta_{\iota_g(h)}$ and $\phi_i^h$ its associated eigenfunction.

We set once again $$F_i : \Cl S^{k+2}\cx (a,b) \cx \Cl O\ra \Cl H^{k}$$
defined $$F_i(\phi, \lambda, h) = \Delta_h\phi - \lambda\phi.$$ We set $$Q_i = F^{-1}(0)\subset \Cl S^k\cx (a,b)\cx \Cl O_i$$ and $\pi : Q_i \ra \Cl O_i$ the restriction to $Q_i$ of the natural projection $\Cl S^k\cx (a,b)\cx \Cl O_i\ra \Cl O_i$.

\begin{Lemm}\label{lem:TanQ}
The set $Q_i$ is a Banach submanifold of $\Cl S^k\cx (a,b)\cx \Cl O_i$ whose tangent space in $(\phi,\lambda,g)$ is :
$$T_{(\phi,\lambda,h)}Q = \left\{(v,\eta,\ddh)\in \Cl H^{k+2}\cx \Bb R\cx T_g\Cl O_g ; \int_M \phi v = 0, (\Delta_h-\lambda)v-\eta\phi+(\bd_h F)\ddh = 0\right\}.$$
Moreover, the projection $\pi$ is an index 0 Fredholm map.
\end{Lemm}

\begin{proof}
By construction, for any $h\in\Cl O_i$, $\Delta_h = \Delta_{\iota_g(h)}$ admits a unique eigenvalue $\lambda_i(h)\in (a,b)$ which is simple : $(\phi, \lambda, h)\in F_i^{-1}(0)$ if and only if $\lambda = \lambda_i(\iota_g(h))$ et $\phi = \phi_i^h$. By Lemma \ref{lem:CaracValSimp},  $0$ is a regular value and by Implicit Functions Theorem, $Q_i$ is a Banach submanifold. By definition, $$T_{(\phi,\lambda,h)}Q = \left\{(v,\eta,\ddh)\in T_\phi\Cl S^{k+2}\cx \Bb R\cx T_g\Cl O_g ;  T_{(\phi,\lambda,g)}F(v,\eta,h) = 0\right\},$$ which gives immediately the expression of the Lemma. As $Q_i$ is a smooth submanifold, $\pi$ is differentiable. We have just seen it is a bijection from $Q$ to $\Cl O_g$, it is therefore a Fredholm operator of index 0.
\end{proof}

To characterize Morse property, we will use the notion of transversality :

\begin{Defi}
Two submanifolds $\Cl N$ of $\Cl N'$ a Banach manifold $\Cl M$ are \ind{transverse} if for any $x\in \Cl N\cap \Cl N'$, $T_x\Cl N$ contains a \emph{closed complement} of $T_x\Cl N'$ in $T_x\Cl M$. Let $\phi : \Cl U\ra \Cl M$ be a $\Cl C^1$ map betwenn two Banach manifold, $\phi$ is transverse to $\Cl N$ if for any point $x\in \Cl U$ such that $\phi(x)\in \Cl N\subset \Cl M$, $(d\phi)^{-1}(T_{\phi(x)}\Cl N)$ admits a \emph{closed complement} in $\Cl U$ and $Im(D_x\phi)$ contains a closed complement of $T_{\phi(x)}\Cl N$ in $T_{\phi(x)} \Cl M$.
\end{Defi}

We note $\xi_0\subset T^*M$ the null section of the cotangent bundle. One can immediately check that a function $f : K\Bb R$ is Morse on $K$ if and only if for any $x\in K, df(x)$ is transverse to $\xi_0$. Therefore, Proposition \ref{prop:KMorseDense} is implied by the following proposition and Baire Theorem :

\begin{Prop}\label{prop:KMorseResiduel}
There is a subset $\Cl M_i\subset \Cl O_i$, generic for the $\Cl C^{k+2}$-uniform topology, such that for any $h\in\Cl M_i$, if $(\phi,\lambda,h)\in\pi^{-1}(h)$, then $d\phi$ is transverse to $\xi_0$ on $K$.
\end{Prop}

\begin{proof}
To prove this result, we will use as in \cite{Uhl76} the following version of Sard-Smale Theorem:

\begin{Theo}[Sard-Smale Theorem 2]\label{th:SardSmale2}
Let $Q,B,X,Y$ et $Y'$ be Banach manifolds, with $Y'\subset Y$, where $X, Y$ and $Y'$ are finite dimensional. Let $\pi : Q\ra B$ be a $\Cl C^k$ Fredholm map with index 0. Then if $f : Q\cx X\ra Y$ is $\Cl C^k$ for $k>\max(1, \dim X+\dim Y'-\dim Y)$, and if $f$ is transverse to $Y'$, then the set
$$\left\{b\in B ; f_b = f|_{\pi^{-1}(b)\cx X} \mbox{ is transverse to } Y'\right\}$$
is generic in $B$.
\end{Theo}
\begin{proof}
This is proven from the original Sard-Smale Theorem in \cite{Uhl76}, p1061.
\end{proof}

\begin{Coro}\label{prop:TransvMorse}
Let us set $\beta :Q\cx K\ra T^*M$ defined by $\beta(\phi,\lambda,h,x) = d\phi(x)$. Assume $\beta$ is $\Cl C^1$ and transverse to the null section $\xi_0$ on $K$. Then the set of $h\in \Cl O_i$ such that $\phi_i^h$ is Morse on $K$ is generic $\Cl O_g$.
\end{Coro}

\begin{proof}
Let $h\in \Cl O_i$, by definition $\phi_i^h$ is Morse on $K$ if and only if $d\phi_0^h$ is transverse to $\xi_0$
on $K$. Now, $\beta_h = \beta|_{\pi^{-1}(h)\cx \inter{K}} = d\phi_i^h$: this corollary is directly implied by the previous version of the Sard-Smale Theorem, with $f=\beta$ and $Y' = \xi_0|_K$. Since the dimension of $M$ as well as the dimension of $\xi_0$ are equal to half the dimension of $T^*M$, it is enough for $\beta$ to be $\Cl C^1$.
\end{proof}

To finish the proof of Proposition \ref{prop:KMorseResiduel}, we have to show that $\beta : Q\cx K \ra T^*M$ is $\Cl C^1$ and transverse to $\xi_0$. As for any $h\in\Cl O_g$, the metric $i_g(h)$ is $\Cl C^{k+2}, k\geq 1$, the Laplacien in local coordinates (see Paragraph \ref{ssec:BasicSpec}) defines an elliptic operator whose coefficients are $\Cl C^{k+1}$. By the local theory of elliptic equations (see \cite{GilTru01}, p185), the eigenfunctions $\phi_i^h$ are $\Cl C^{k+2}$. Therefore, the $\beta$ is $\Cl C^{k+1}$. Let us show now it is transverse to $\xi_0$. 

By definition, the image of $(\phi,\lambda,h,x)\in Q\cx K$ by $\beta$ is in $\xi_0$ if and only if $d\phi(x) = 0$. 
Let $(\phi,\lambda,h,x)\in Q\cx K$ be such that $d\phi(x) = 0$, we want to show that
$$Im(d_{(\phi,\lambda,h,x)}\beta)+ T_{(x,d\phi(x))}\xi_0 = T_{(x,0)}(T^*M).$$

Let $(v,\eta,\ddh)\in T_{(\phi,\lambda,h)}Q$. For any $X\in T_xM$, we have
 $$d\beta(v,\eta,\ddh,X) = d_x v + d_x^2\phi(X,.),$$
where $d_x^2\phi = \nabla d\phi$ is the second differential of $\phi$ in any local coordinates basis, as $x$ is a critical point of $\phi$.
Let $U_x$ be a neighbourhood of $x$ (with compact closure) diffeomorphic to $\Bb R^n$ and $(x_1,...,x_n)$ be local coordinates on this neighbourhood. We only have to show that for any $i=1...n$, there exists $(v,\eta,\ddh)\in T_{(\phi,\lambda,g)}Q$ such that 
$\bd_{x_i}v\neq 0.$

By Lemma \ref{lem:TanQ}, the tangent space $T_{(\phi,\lambda,h)}Q$ is characterized by the following equalities :
$$\int_M \phi v dv_g = 0\hbox{ et } (\Delta_g-\lambda)v-\eta\phi+(\bd_h F)\ddh = 0.$$
If we restrict to variations on $Q$ which preserve $\lambda_i$, that is $\eta = 0$, these constraints become :
$$(\Delta_h-\lambda)v = -(\bd_h F)\ddh,$$
with $v$ orthogonal to $\phi$. We will need to express such a $v$ from a \emph{Green function associated to $\phi$}.

\begin{Prop}\label{prop:GreenAdapt}
Let $(\phi,\lambda,h)\in Q_i$. There a $\Cl C^2$-function $G_\phi : M\cx M\bs \mbox{Diag}\ra \Bb R$ on $M\cx M\bs \mbox{Diag}$ satisfying the following properties :
\bi
\item for any $x\neq y$, $G_\phi(x,y) = G_\phi(y,x)$ ;
\item for any $x\in M$, $$\int_M\phi(y)G_\phi(x,y)dv_h(y) = 0 $$
and $G_\phi\in \Cl H^1(M)$ ;
\item for any $x\in M$ and $f\in\Cl L^2(M)$ such that $\int_M \phi f dv_h = 0$, let us set $u(x) = \int_M G_\phi(x,y)f(y) dv_h(y)$, we have :
\bi
\item $$\nabla_h u(x) = \int_M \nabla_h G_\phi(x,y)f(y) dv_h(y) ;$$
\item $$(\Delta_h-\lambda)u(x) = f(x).$$
\ei
\ei
We have noted $\nabla_h$ the connexion associated to $\iota_g(h), \Delta_h$ its Laplacian and $dv_h$ its volume element on $M$. We call $G_\phi$ the \ind{modified Green function} associated to $\phi$.
\end{Prop}

\begin{proof}
Let $(\phi,\lambda,h)\in Q_i$. As $\iota_g(h)$ has bounded geometry at order $k\geq 1$, by the Theorem 4.2 of \cite{Dod83}, $(M,\iota_g(h))$ admits a unique $\Cl C^k$ \ind{heat kernel} $p : M\cx M\cx ]0,\infty[\ra \Bb R^*_+$, that is the fundamental solution of the equation
$$\Delta f = -\frac{\bd f}{\bd t}.$$ Recall that we assumed that $\iota_g(h)$ has $q+1$ first distinct eigenvalues
$$0<\lambda_0<...<\lambda_i = \lambda<...<\lambda_q<\lambda_0^{ess}(\iota_g(h)) = \lambda_0^{ess}(g).$$
For all $m\in [0,q]$, we note $\phi_m$ the eigenfunction associated to $\lambda_m$, by construction $\phi = \phi_i$. The Spectral Decomposition Theorem (see \cite{ReeSim80}, chapter VIII) for the heat operator $e^{-t\Delta}$ implies that for any $(x,y,t)\in M\cx M\cx \Bb R^*_+$,
$$p(x,y,t)= \sum_{m = 0}^q e^{-\lambda_m t}\phi_m(x)\phi_m(y) + p_{q+1}(x,y,t),$$
whith the following properties :
\bi
\item for any $(y,t)\in M\cx \Bb R^*_+$ and any $m\in [0,q], p_{q+1}(.,y,t)$, is $L^2$-orthogonal to $\phi_m$ ;
\item let $\lambda_{q+1}$ be the $q+2$-th discrete eigenvalue of $\Delta_h$ if it exists, otherwise $\lambda_{q+1} = \lambda_0^{ess}(g)>\lambda_i(\iota_g(h))$, then for any $(x,y)\in M\cx M$,
$$p_{q+1}(x,y,t) = \Cl O(e^{-\lambda_{q+1}t})$$
when $t\ra\infty.$
\ei
We set
$$p_{i+1}(x,y,t) = p(x,y,t)-\sum_{m = 0}^i e^{-\lambda_m t}\phi_m(x)\phi_m(y).$$

\begin{Lemm}
The function $e^{\lambda t} p_{i+1} : M\cx M\cx ]0,\infty[\ra \Bb R$  defined for all $(x,y,t)\in M\cx M\cx ]0,\infty[$ by
$$e^{\lambda t} p_{i+1} = e^{\lambda t}p_{i+1}(x,y,t) = e^{\lambda_i t}p(x,y,t)-\sum_{m = 0}^i e^{(\lambda_i-\lambda_m) t}\phi_m(x)\phi_m(y)$$
satisfies the following properties :
\bi
\item $e^{\lambda t} p_{i+1}$ is solution of the following modified heat equation : for any $(x,y,t)\in M\cx ]0,\infty[,$ we have
$$(\Delta_h - \lambda)p_\phi (x,y,t) = - \frac{\bd}{\bd t} p_\phi(x,y,t),$$
where the Laplacian is taken in the first space coordinate ;

\item for any $(x,t)\in M\cx ]0,\infty[,$ $e^{\lambda t} p_{i+1}(x,.,t)$ is orthogonal to $\phi$ :
$$\int_M e^{\lambda t} p_{i+1}\phi(y)dv_g(y) = 0 ;$$

\item for any bounded summable function $f : M\ra \Bb R$ such that $\int_M f\phi = 0$,
the function
$$x\fa\int_Me^{\lambda t} p_{i+1}(x,y,t)f(y)dv_g(y)$$
satisfies the same modified heat equation, and we have
$$\lim_{t\ra 0}\int_Me^{\lambda t} p_{i+1}f(y)dv_g(y) = f(x) ;$$
\item for any $(x,y)\in M\cx M$, we have
$$p_\phi(x,y,t)\sim C(x,y) e^{-\eta t}$$
when $t\ra\infty$, where $\eta = \lambda_{i+1}-\lambda>0$.
\ei
\end{Lemm}

This lemma comes straightforward from the spectral decomposition of the heat kernel, and the classical fact that for any $x\in M$,
$$\lim_{t\ra 0}p(x,y,t) = \delta_x,$$
where $\delta_x$ is the Dirac mass in $x$ (see \cite{Cha84}, chapter VIII).

For any $x,y\in M, x\neq y$, we set now
$$G_\phi(x,y) = \int_0^\infty e^{\lambda t} p_{i+1}(x,y,t)dt.$$
By the previous lemma, $G_\phi$ is well defined on $M\cx M \bs Diag$ and is $\Cl C^{k+1}$ on $M\cx M \bs Diag$. The symetry of $G_\phi$ and its orthogonality to $\phi$ comes from corresponding properties of $e^{\lambda t} p_{i+1}$. It is a classical result (see for example Lemma 3.2 of \cite{Dod83}) that, when $y\ra x$, we have the following equivalent :
$$G_\phi(x,y)\sim \left\{\begin{array}{cc}
														\frac{C(x)}{d_g(x,y)^{n-2}} & \mbox{ if }n>2,\\
														C(x)\ln(d_g(x,y)) & \mbox{ if }n=2\end{array}\right.$$
where $d_g(x,y)$ is the distance between $x$ and $y$ induced by $g$. Using this equivalent and the fact that for any $t>\epsilon$ and $x\in M$,
\beq\label{eq:L2Pphi}
\int_M e^{\lambda_i t} p_{i+1}(x,y,t)^2 dv_g(y)\leq e^{-\eta(t-\epsilon)}\int_M e^{\lambda t} p_{i+1}(x,y,\epsilon)dv_g(y),\eeq computing separately the integral on a neighbourhood of $(x,0)$ and on the rest of $M\cx [\epsilon,\infty[$, we show then that for any $f\in \Cl L^2(M)$ orthogonal to $\phi$, the integral $$u(x) = \int_M G_\phi(x,y)f(y) dv_g(y)$$ is well defined, and is the unique solution of the equation $(\Delta-\lambda)u = f$ orthogonal to $\phi$. An analogous estimates to (\ref{eq:L2Pphi}) for the $L^2$-norm of the gradient of  $e^{\lambda_i t} p_{i+1}$ shows the gradient of $u$ is given by the integral form in point 3.(b) and that $G_\phi(.,y)\in \Cl H^1(M)$.
\end{proof}

We go back to prove that $\beta$ is transverse to $\xi_0$. Recall we note $\bd_h F$ the derivative of $F$ along $\Cl O_i$ an $$J = Im(\bd_h F)\subset \Cl H^{k}(M).$$ As $\Delta_h-\lambda$ is invertible from the orthogonal of $\phi$ into itself ($\lambda = \lambda_i$ is a simple eigenvalue), for any $j= (\bd_h F)\ddh\in J\cap \{\phi\}^\bot$, there exists a unique $v\in \{\phi\}^\bot$ such that $(v,0,\ddh)\in T_{(\phi,\lambda,h)}Q$ and 
$$(\Delta_h-\lambda)v = -(\bd_h F)\ddh.$$ 
This $v$ is explicitely given by
$$v(x) = -\int_MG_\phi(x,y)j(y)dv_h(y),$$ where $G_\phi$ is the modified Green function associated to $\phi$ constructed in the previous proposition. We saw then that
$$\bd_{x_i}v(x) = -\int_M(\bd_{x_i}G_\phi)(x,y)j(y)dv_h(y).$$

Assume now by contradiction that for any $(v,0,\ddh)\in T_{(\phi,\lambda,h)}Q$, $$\bd_{x_i}v(x) = 0.$$ 

By the previous argument, for any $j\in J$ orthogonal to $\phi$, 
$$\int_M(\bd_{x_i}G_\phi)(x,y)j(y)dv_h(y) = 0.$$
Set $$a = \frac{\int_M(\bd_{x_i}G_\phi)(x,y)\phi(y) dv_h(y)}{\int_M\phi^2dv_h},$$
we have then for any $j\in J$,
$$\int_M\left[(\bd_{x_i}G_\phi)(x,y)-a\phi(y)\right]j(y)dv_h = 0.$$
But we proved in the Lemma \ref{lemm:JDense} that this implies that for any $y\in K$,
$$\bd_{x_i}G_\phi(x,y) = a\phi(y) + cste.$$ This is impossible for, by Proposition \ref{prop:GreenAdapt}, $\bd_{x_i}G_\phi(x,.)$ is singular in $x$. This concludes the proof of Proposition \ref{prop:KMorseResiduel}, which itself implies the Proposition \ref{prop:KMorseDense}.
\end{proof}
\end{proof}

Propositions \ref{prop:KMorseOuvert} and \ref{prop:KMorseDense} show that for any $i\in[0,q]$ and any compact $K\subset M$, the set $\Cl M_i(K)\subset \Cl U$ of the metrics in $\Cl U$ such that for any $g\in\Cl M_i(K)$, the $i$-th eigenfunction $\phi_i^g$ is Morse on $K$ \emph{is open and dense in $\Cl U$} for the $\Cl C^{k+2}$-strong topology. As $M$ is a countable union of compact sets, this ends the proof of Theorem \ref{th:PhiM}.
\end{proof}

\section{Fundamental domains and bottom of the spectrum}

Let $(M,g)$ be a riemannian manifold of dimension $n$. We say closed submanifold with boundary $D\subset M$ of dimension $n$ has \ind{piecewise $\Cl C^1$ boundary} if there exists a $\Cl C^1$-atlas of $M$ into $\Bb R^n$ such that the image of any open set of $D$ is open in the quadrant $$\{(x_1,...,x_n,) ; x_1\geq 0,...,x_n\geq 0\}.$$ The boundary $\bd D$ itself is then a locally finite union of $n-1$-dimensional closed submanifolds with piecewise $\Cl C^1$ boundary.

Let $p : M'\ra M$ be a riemannian covering. We will call a \ind{fundamental domain} for $p$ a \emph{closed connected submanifold $D\subset M'$ with piecewise $\Cl C^1$ boundary} such that $p(D) = M$ and each element of $M\bs p(\bd D)$ has a unique preimage in $D$ by $p$. 

\subsection{Morse functions and fundamental domains}

Let $p : M'\ra M$ be a riemannian covering, to any Morse function on $M$ one can associate a fundamental domain which is adapted to its \ind{gradient flow} by the following theorem.

\begin{Theo}\label{theo:DomFondMorse}
Let $p : M'\ra M$ be a riemannian covering and $f : M\ra \Bb R_+^*$ a Morse function such that, for any $a>0$, 
$$M_a = \left\{x\in M ; f(x)\geq a\right\}$$
is compact. Let us note $\tilde{f} = f\circ p$, there exists a fundamental domain $D\subset M'$ for $p$ such that $\nabla \tilde{f}$ is tangent to $\bd D$. Moreover, if $f$ has a finite number of critical points, then for any $1\leq k\leq n,$ the boundary $\bd D$ has a finite number of $k$-dimensional smooth components.
\end{Theo}

This is a special case of Theorem \ref{theo:DomFondStrat}, which we prove just below. The fundamental domain constructed by this result is stable under the gradient flow of $\tilde{f}$. In other words, $f$ satisfies Neumann boundary conditions on $\bd D$. 

\begin{Defi}\label{def:GradStrat}
Let $(M,g)$ be a riemannian manifold and $f : M\ra \Bb R$ a $\Cl C^1$ function. We say $f$ has a \ind{simply stratified gradient} if $M$ can be decomposed into
$$M = \coprod_i M_i,$$
where the $M_i$ are \emph{disjoint open $\Cl C^1$ submanifolds}, \emph{simply connected} and \emph{stable by the gradient flow} of $f$  (in the metric $g$), and the partition $(M_i)_i$ is \emph{locally finite}.
\end{Defi}

A locally finite partition of a topological set into manifolds is called a \ind{stratification} (see for instance \cite{GorMac88} p. 37), which gives the name of this property. The elements $M_i$ of the partition are called \ind{strata} ; the closure of a strata is a submanifold with piecewise $\Cl C^1$ boundary.

\begin{Theo}\label{theo:DomFondStrat}
Let $p : M'\ra M$ be a riemannian covering, $f : M\ra \Bb R$ a function with simply stratified gradient and $\tilde{f} = \pi\circ f$. There exists a fundamental $D\subset M'$ for $p$ such that $\nabla \tilde{f}$ is tangent to $\bd D$.
\end{Theo}

\begin{proof}
Let $p : M'\ra M$ be a riemannian covering and $f : M\ra \Bb R$ a function with simply stratified gradient and $M = \coprod_i M_i$ the associated stratification. We will note $f := \tilde{f}$. Let $M_1\subset M$ be a stratum with codimension $0$, as $M_1$ is simply connected, $p$ is trivial over $M_1$ : all connected components of $p^{-1}(M_1)$ are diffeomorphic to $M_1$. We call each of them a \ind{lift} of $M_1$ by $p$. Let $M'_1$ be one of these lifts. The boundary of $N'_1$ is a locally finite union of strata with codimension at least $1$. If $\bar{M_1} = M$, we set $D = \bar{M'_1}$ : one can immediately check that it is a fundamental domain for $p$ such that $\nabla f$ is tangent to $\bd D$. If $\bar{M_1}\neq M$, there exists $M_2\neq M_1$ another stratum with codimension 0 such that $\bar{M_1}$ and $\bar M_2$ contain a common codimension 1 stratum $M_{12}$. This latter is completely included in the interior of $\bar{N_1\cup N_2}$, which is therefore connected. Let $M'_2$ be a lift of $M_2$ such that the interior of $\bar{M'_1\cup M'_2}$ contains a lift of $N_{12}$. If $\bar{M_1\cup M_2} = M$, then $D=\bar{M'_1\cup M'_2}$ is a fundamental domain satisfying the conclusions of the theorem. Otherwise, there exists a codimension 0 stratum $M_3$ disjoint from $M_1$ et $M_2$ such that $\bar{M_1\cup M_2}$ and $\bar{M_3}$ contain a common codimension 1 stratum $M_{23}$. We repeat the previous argument, and as there is at most a countable number of codimension 0 strata, repeting the argument a countable number of time we obtain a fundamental domain satifying the properties we wanted.
\end{proof}

The Theorem \ref{theo:DomFondMorse} is a corollary of the previous theorem and the following result:

\begin{Theo}[Thom]
Let $(M,g)$ be a riemannian manifold and $f : M\ra \Bb R_+^*$ a Morse function such that for any $a>0$, 
$$M_a = \left\{x\in M ; f(x)\geq a\right\}$$
is compact. Then $f$ has a simply stratified gradient, with a finite number of strata if $f$ has a finite number of critical points.
\end{Theo}

\begin{proof}
Let $(M,g)$ be a riemannian manifold and $f$ such a Morse function. All critical points of $f$ are isolated. Recall the \ind{stable manifold} of a critical point $\alpha$ is the set of points $x\in M$ such that the gradient line of $f$ through $x$ (oriented towards increasing $f$) ends in $\alpha$. It can be shown that for any critical point $\alpha$ of $f$ with index $r$, the stable manifold of $\alpha$ is diffeomorphic to $\Bb R^{n-r}$ (this result, due to Thom, is proved for example in \cite{AbRob67}, p87). Moreover, by assumption for any $a>0$, 
$M_a = \left\{x\in M ; f(x)\geq a\right\}$
is compact. This implies that any point belong to a (unique) stable manifold, and that the set of these stable manifolds is locally finite. Stable manifolds give then a stratification of $M$ into simply connected submanifolds, which are by definition stable by the gradient flow of $f$. 
\end{proof}


\begin{center}\begin{figure}\label{fig:DomMorse}
\includegraphics[width = \textwidth]{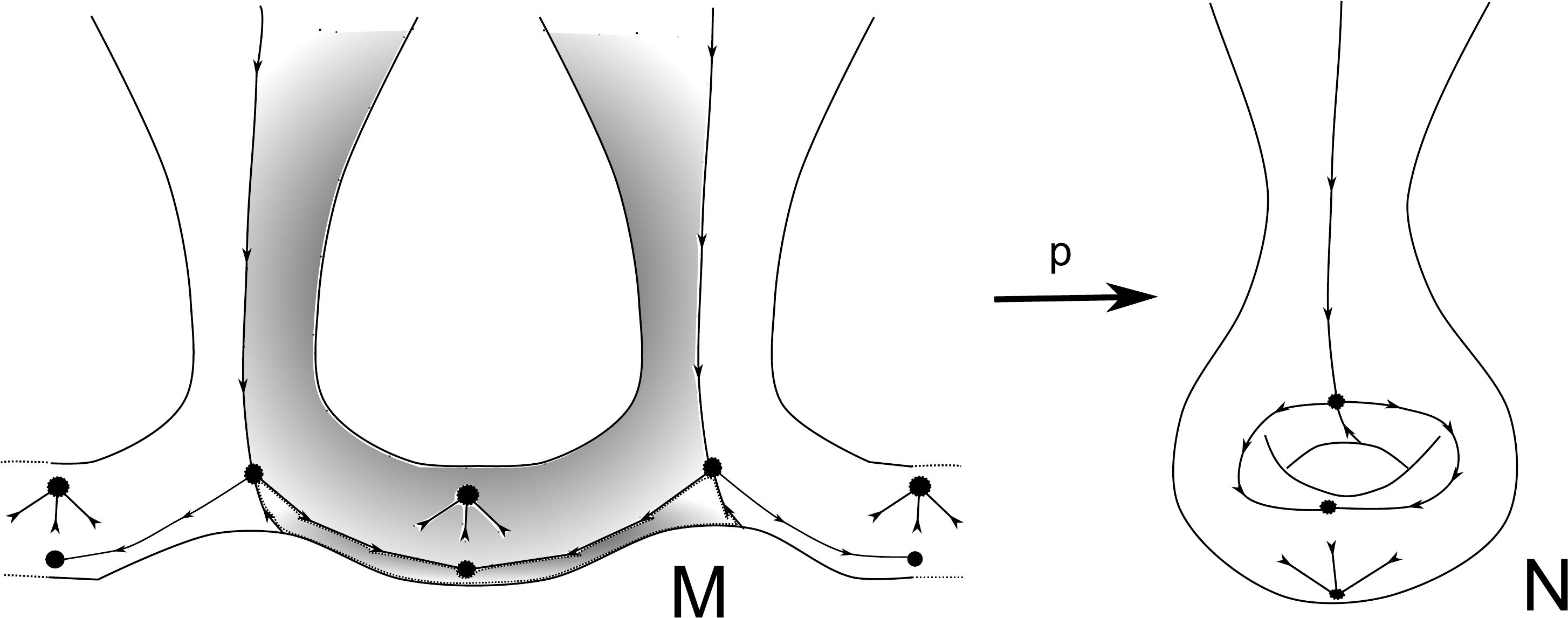}
\caption{\emph{Stable manifolds and fundamental domain :} we have drawn the critical points of the function and their stable manifolds. The fundamental domain correspondinf to the stable manifold of a local maximum is shaded.}
\end{figure}\end{center}

A Morse function satisying the hypothesis of Theorem \ref{theo:DomFondStrat} always exists, by a classical result of Morse Theory :

\begin{Theo}
On any differentiable manifold $M$ there exists a Morse function $f : M\ra \Bb R^*_+$ such that for any $a>0$, 
$$M_a = \left\{x\in M ; f(x)\geq a\right\}$$
is compact.
\end{Theo}

This is Corollary 6.7 of \cite{Mil63}, whose proof is detailed all along chapter 6 of J. Milnor's book.

\subsection{Bottom of the spectrum through fundamental domains}

Let $(M,g)$ be a riemannian manifold with bounded geometry at order $k\geq\frac{n}{2}$ and infinite volume. We assume the bottom of the spectrum of $\Delta_g$ is an isolated eigenvalue, whose (normalized) eigenfunction we call $\phi_0$. Let $p : M'\ra M$ be a riemannian covering. We now want to apply the genericity theorem of the first section and the construction of fundamental domains associated to Morse functions we have just presented to prove Theorem \ref{th:SpecDom} :

\begin{Theo}\label{th:DomFond}
With previous notations, 
$$\lambda_0(M) = \sup_D \lambda_0(D),$$
where $D\subset M'$ runs over all fundamental domains for $p$. Any fundamental domain $D$ on which the lift of $\phi_0$ satisfies Neumann boundary conditions satisfies $\lambda_0(D) = \lambda_0(M)$.
\end{Theo}

Recall $\lambda_0(D)$ is the bottom of the spectrum of $D$ avec condition de Neumann au bord. A first half of our theorem is immediate :

\begin{Lemm}\label{lem:MinL}
For any fundamental domain $D\subset M'$  for $p$, $\lambda_0(M)\geq \lambda_0(D).$
\end{Lemm}

\begin{proof}
By the Min-Max principle, $\lambda_0(M)$ is the infimum of the Rayleigh quotients of $\Cl H^1$ functions on $M$, and for Neumann boundary conditions, $\lambda_0(D)$ is the infimum of the Rayleigh quotients of $\Cl H^1$ functions on $D$ without any boundary restriction. Each $\Cl H^1$ function on $M$ can be lifted to a $\Cl H^1$ function on $D$, which gives the lemma.
\end{proof}

To apply Theorem \ref{theo:DomFondStrat} to construt a domain adapted to $\phi_0$, we need to control its level sets. 

\begin{Lemm}\label{lem:PhiCompact}
For any $a>0$, the set $M_a = \left\{x\in M ; \phi_0(x)\geq a\right\}$ is compact.
\end{Lemm}

\begin{proof}
Let $a>0$, we assume by contradiction $M_a$ is not compact. As $M$ is complete and $\phi_0$ is continuous, there exists a sequence  $(x_n)_{n\geq 0}\in (M_a)^{\Bb N}$ of points in $M_a$, and $r>0$ such that the $B(x_n,r)$ are disjoint to each other. We assume $r$ to be less than the injectivity radius $inj(M)$, which we assumed to be positive. By definition of $M_a$, for any $n\in \Bb N, \phi_0(x_n) \geq a$. Moreover, as the Ricci curvature is also bounded by assumption, by \cite{Yau75} there exists a constant $C$ which only depends on the dimension of $M$ and on a lower bound on its Ricci curvature such that the \emph{positive} solutions of $\Delta f = \lambda f$ satisfy $$\left|\frac{\nabla f}{f}\right|\leq C.$$
Therefore, there exists $r_0<r$ such that for any $n\geq 0$ and any $x\in B(x_n,r_0)$,
$$\phi_0(x)\geq \frac{a}{2}.$$
Moreover, it is shown in \cite{GrePet92} that when $r_0$ is less than the injectivity radius of $M$, there exists a constant $C_d$ which only depends on the dimension $d$ of $M$ such that for any $n\geq 0$, $$\Vol(B(x_n,r_0))>C_d r_0^d>0.$$ Therefore,
$$\int_N\phi_0^2 \geq \sum_{n\geq 0} \frac{a^2}{4}\Vol(B(x_n,r_0)) = \infty,$$
which is a contradiction as $\phi_0$ is $L^2$.
\end{proof}


This result will allow us to prove Theorem \ref{th:SpecDom} when the metric is generic :

\begin{Prop}\label{prop:DomOptiMorse}
If $\phi_0$ is Morse on $M$, then there exists a fundamental domain $D\subset M$ with piecewise $\Cl C^1$ boundary such that the lift of $\phi_0$ satisfies Neumann boundary conditions on $\bd D$ and 
$\lambda_0(D) = \lambda_0(M).$
\end{Prop}

\begin{proof}
From Theorem \ref{th:SullNeum}, if there exists $D$ such that $\phi_0$ satisfies Neumann boundary conditions on $\bd D$, then as $\phi_0$ is positive, $\lambda_0(D)\geq\lambda_0(M)$, which ensure $\lambda_0(D) = \lambda_0(M)$ by Lemma \ref{lem:MinL}. Since the sets $M_a$ are compact for any $a>0$, such a fundamental domain exists as soon as $\phi_0$ is Morse by Theorem \ref{theo:DomFondMorse}. 
\end{proof}

Note that the proof of this proposition does not require that the geometry of $(M,g)$ be bounded at order $k\geq \frac{n}{2}$, but only at order 0. When $\phi_0$ is not Morse, we prove Theorem \ref{th:SpecDom} by approximation, which will use the complete boundedness assumption through the genericity Theorem \ref{th:PhiM} :

\begin{Prop}\label{prop:DomPresqueBon}
For any $\epsilon>0$, there exists a fundamental domain $D_\epsilon$ with piecewise $\Cl C^1$ boundary such that
$$\lambda_0(D_\epsilon)\geq \lambda_0(M) - \epsilon.$$
\end{Prop}

\begin{proof}
Let $\epsilon>0$, by Proposition \ref{prop:SurharmPos}, it is enough to show there exists a fundamental domain $D_\epsilon$ which admits a $(\lambda_0(M,g_0)-\epsilon)$-superharmonic positive function for the Laplacian associated to the metric $g$, and which satisfies Neumann boundary conditions on $\bd D_\epsilon$. We still note $\lambda_0 = \lambda_0(g)$ et $\phi_0$ the associated eigenfunction of $\Delta_{g}$. 

\begin{Lemm}\label{lem:PresqHarm}
For any $\epsilon>0$, there exists a neighbourhood $\Cl V_\epsilon$ of $g$ for the $\Cl C^2$-strong topology such that for any $h\in\Cl V_\epsilon$, the eigenfunction $\phi_0^h$ of $\Delta_h$ is $(\lambda_0-\epsilon)$-superharmonic positive for the Laplacian $\Delta_g$.
\end{Lemm}

\begin{proof}
By Theorem \ref{th:StabTrou}, there exists a neighbourhood $\Cl V_1$ dofe $g$ such that for any $h\in\Cl V_1$, the bottom of the spectrum of $\Delta_h$ is a positive isolated eigenvalue and $$|\lambda_0(g_0)-\lambda_0(g)|\leq \epsilon.$$

We will use local coordinates which are adapted to the computation of the Laplacian for the metric $g$ : the \emph{harmonic coordinates.}
\begin{Defi}
Let $(x_i)_{i = 1,...,n}$ be a local coordinates basis. It is said to be \ind{harmonic coordinates} for the metric $g$ if for any $i = 1,...,n$, $$\Delta_{g}(x_i) = 0.$$
\end{Defi}

Such coordinates always exists on  a small open set. The following theorem will allow us to use them to define the $\Cl C^2$-strong topology :

\begin{Theo}
Let $(M,g)$ be a riemannian manifold with bounded geometry at order 1. Then there exists a locally finite atlas whose charts all define harmonic coordinates.
\end{Theo}
This is a straightforward consequence of Theorem 1.3 of \cite{Heb96}. 

We consider from now on a fixed harmonic atlas $\g U$ for the metric $g$. In any local chart $(x_i)$, for any $\Cl C^2$ function $f$, the Laplacian can be written
$$\Delta_h f = (h^{ij}(\bd_{ij} f - \Gamma_{ij}^p(h)\bd_p f),$$
where the $\Gamma_{ij}^p(h)$ are the Christoffel symbol of $h$ in the local chart $(x_i)$, and where we used Einstein summation convention (see \cite{Heb96}, Chapter 1). We have in particular for any $p = 1,...,n$,
$$\Delta_{g} x_p = g^{ij}\Gamma_{ij}^p(g),$$
still using Einstein summation convention. In a local basis of harmonic coordinates for $g$, we have then for any $\Cl C^2$ function $f$ 
$$\Delta_{g} f = g^{ij}\bd_{ij} f,$$
which implies
$$(\Delta_h-\Delta_g)f = (h^{ij}-g^{ij})\bd_{ij}f  - h^{ij}\Gamma_{ij}^p(h)\bd_p f.$$
For $f = \phi_0^h$, we have then
$$\Delta_h \phi_0^h = h^{ij}(\bd_{ij} \phi_0^h - \Gamma_{ij}^p(h)\bd_p \phi_0^h) = \lambda_0(h) \phi_0^h$$
therefore
\beq\label{eq:MajEps1}
|h^{ij}\bd_{ij}\phi_0^h - \lambda_0(h)\phi_0^h|\leq |h^{ij}\Gamma_{ij}^p\bd_p \phi_0^h|.\eeq

By continuity of the application which, in a fixed coordinate basis, maps the metric to its Christoffel symbols, there exists a neighbourhood of $g$ for the $\Cl C^2$-strong topology (defined through the harmonic atlas harmonique we have fixed) $\Cl V_2\subset V_1$ such that for any $h\in\Cl V_2$, in any local chart of $\g U$ and for any $p = 1,...,n$, we get
$$|h^{ij}\Gamma_{ij}^p(h)|\leq \epsilon.$$

Moreover, by assumption for any $h\in\Cl V_1$, the Ricci curvature of $(M,h)$ and its injectivity radius are bounded by uniform constants on $\Cl V_1$. As $\phi_0^h$ is a positive eigenfunction on $\Delta_h$, by \cite{Yau75}, there exists a constant $C$ depending only on $\Cl V_1$ (and not on $h$) such that for any $h\in\Cl V_1$, we get
$$|\nabla\phi_0^h|\leq C\phi_0^h.$$
The upper bound (\ref{eq:MajEps1}) becomes then 
$$|h^{ij}\bd_{ij}\phi_0^h - \lambda_0^h\phi_0^h|\leq C\epsilon \phi_0^h,$$
therefore as $\Cl V_2\subset \Cl V_1$,
$$|h^{ij}\bd_{ij}\phi_0^h - \lambda_0\phi_0^h|\leq (C+1)\epsilon \phi_0^h.$$
This becomes then
$$|\Delta_{g}\phi_0^h - \lambda_0\phi_0^h| = |g^{ij}\bd_{ij}\phi_0^h - \lambda_0\phi_0^h|\leq |(g^{ij}-h^{ij})\bd_{ij}\phi_0^h|+|h^{ij}\bd_{ij} - \lambda_0\phi_0^h|.$$
Moreover, there exists a neighbourhood $\Cl V_3\subset V_2$ of $g$ for any $h\in \Cl V_3$, in any chart of $\g U$ and for all $i,j = 1,...,n$, we get
$$|h^{ij}-g^{ij}|\leq \epsilon |h^{ij}|.$$
We have then 
$$|(g^{ij}-h^{ij})\bd_{ij}\phi_0^h|\leq \epsilon|h^{ij}\bd_{ij}\phi_0^h| = \epsilon|\lambda_0^h \phi_0^h - h^{ij}\Gamma_{ij}^p\bd_ p\phi_0^h|,$$
therefore $$|(g^{ij}-h^{ij})\bd_{ij}\phi_0^h|\leq \epsilon(\lambda_0+\epsilon+C\epsilon)\phi_0^h$$
by previous arguments. Eventually, for any metric $h\in\Cl V_3$,

$$|\Delta_{g}\phi_0^h - \lambda_0\phi_0^h|\leq (C+1+\lambda_0+C\epsilon)\epsilon \phi_0^h.$$
Up to replacing $\epsilon$ by $B\epsilon$ where $B$ is a positive constant which do not depend on $\epsilon$, this ends the proof of Lemma \ref{lem:PresqHarm}.
\end{proof}

Let $\Cl V_\epsilon$ be the neighbourhood of $g_0$ for the $\Cl C^2$-strong topology given by the previous lemma. For any $k\geq 0$, it is also a neighbourhood of $g$ for the $\Cl C^{k+2}$-strong topology. By Theorem \ref{th:PhiM}, there exists $h\in\Cl V_\epsilon$ such that $\phi_0^h$ is Morse on $M$. By Proposition \ref{prop:DomOptiMorse}, there exists then a fundamental domain $D_\epsilon$ for $p$ such that $\phi_0^h$ satisfies Neumann boundary conditions on $\bd D_\epsilon$. By the previous lemma, the fonction $\phi_0^h$ is $\lambda_0-\epsilon$-superharmonic for the Laplacian $\Delta_g$. This ends the proos of Proposition \ref{prop:DomPresqueBon}, which also concludes the proof of Theorem \ref{th:SpecDom} in the general case.
\end{proof}

\subsection{Application : bottom of the spectrum of coverings}\label{sec:Persp}
The study of the bottom of the spectrum of fundamental domains for a covering can be used to study the bottom of the spectrum of the covering itself. We have shown in \cite{Tap09a} the following result :

\begin{Theo}
Let $(M,g)$ be a riemannian manifold whose Laplacian has an isolated first eigenvalue. Let $p : M'\ra M$ be a riemannian normal covering with deck-group $\Gamma$, such that there exists a fundamental domain $D\subset M'$ for $p$ satisfying $\lambda_0(D) = \lambda_0(N)$. Then the bottom of the spectrum of $M'$ satisfies
$$\lambda_0(M')\geq \lambda_0(M)$$
with equality if and only if $\Gamma$ is amenable.
\end{Theo}

By Proposition \ref{prop:DomOptiMorse}, if $(M,g)$ has a positive injectivity radius, a lower bound on its Ricci tensor, and if the first eigenfunction is Morse, then there exists a fundamental domain $D\subset M'$ for $p$ satisfying $\lambda_0(D) = \lambda_0(N)$. From our genericity Theorem \ref{th:PhiM}, we deduce then the following:

\begin{Coro}\label{coro:SpecRevet}
Let $(M,g)$ be a riemannian $n$-manifold with infinite volume and bounded geometry at order $k\geq\frac{n}{2}$. Assume the Laplacian $\Delta_g$ has an isolated first eigenvalue. If $g$ is generic for the $\Cl C^{k+2}$-strong topology, then for any normal covering $p : M'\ra M$ with deck-group $\Gamma$, we have
$$\lambda_0(M) \geq \lambda_0(N)$$
with equality if and only if $\Gamma$ is amenable.
\end{Coro}

This corollary is a generalisation of a result of R. Brooks (cf \cite{Broo85Reine}), who shows that under a technical hypothesis on the fundamental domains of $p : M'\ra M$, which implies among other things that the first eigenvalue of $\Delta_g$ is isolated, then  $\lambda_0(M') \geq \lambda_0(M)$ with equality if and only if the deck-group $\Gamma$ is amenable. See \cite{Tap09a} for more details on these works of Brooks and the connections with our result.
\pgh
Let us note that we only used the genericity of the metric to show that the first eigenfunction has a simply stratified gradient (cf Definition \ref{def:GradStrat}). If we were able to show that for smooth metrics with isolated first eigenvalue, the first eigenfunction has always a simply stratified gradient, then the equality would always be attained in Theorem \ref{th:SpecDom} and the conclusion of Corollary \ref{coro:SpecRevet} would hold for any smooth metric. 

\appendix
\section{Bottom of the spectrum with Neumann boundary conditions}\label{sec:SullNeum}

Let $(M,g)$ be a smooth complete non compact manifold whith piecewise $\Cl C^1$ boundary. We now prove the following characterization of the bottom of the spectrum with Neumann condition, which we used in the proof of Theorem \ref{th:SpecDom} : 

\begin{theo}\label{th:SullNeum}
For any real number $\lambda$, there exists a smooth function $\phi$ which is \emph{$\lambda$-harmonic positive} on $M$ with Neumann boundary conditions on $\bd M$ if and only if $\lambda\leq\lambda_0(M)$, where $\lambda_0(M)$ is the bottom of the spectrum of the Laplacian on $M$ with Neumann boundary conditions.
\end{theo}

We can also state this theorem as in our introduction :

$$\lambda_0(M) = \sup\left\{\lambda\in\Bb R : \exists f\in\Cl C^\infty(C), f>0, \frac{\bd f}{\bd\nu}\equiv 0, \Delta f = \lambda f\right\}.$$

When $\bd M = \vd$, this is exactly Theorem 2.1 of \cite{Sul87} ; it was also proved in \cite{CheYau75}, p 345. We adapt the proof to the case of Neumann boundary condition, using as in \cite{Sul87} Section 3-4, probabilistic tools from the Diffusion Theory which we will try to render readable for geometers. 

\subsection{Brownian motion with reflexion on $\bd M$}

Let $(M_j)_{j\in\Bb N}$ be an \emph{increasing} family of open sets of $M$ with compact closure such that $\bigcup_j M_j = M$, all containing a fixed open set $K\subset \bd M$.
For any $j\geq 0$, we note $\bd^1 M_j = \bd M\cap M_j$ and $\bd^2 M_j = \bd M_j\bs\bd M$ which we also assume to be piecewise $\Cl C^1$.
We will note $$\lambda_0^j=\lambda_0^{N,D}(M_j) = \inf_f\frac{\norm{\nabla f}^2_{L^2(M_j)}}{\norm{f}^2_{L^2(M_j)}}$$
where $f$ runs over all smooth functions with compact support in $M_j=\bar{M}_j\bs \bd^2 M_j$ : $\lambda_0^j$ is the bottom of the spectrum on $M_j$ with \emph{Neumann boundary conditions} on $\bd_1M_j$ and \emph{Dirichlet boundary conditions} on $\bd^2 M_j$.
We have then $\lambda_0(M) = \inf_j\lambda_0^j.$ As the $(M_j)_j$ are compact, from the classical theory of partial differential equations, for any $j\in\Bb N$ there exists an eigenfunction $\phi_0^j\in\Cl H^1(\bar{M_j})$ which satisfies Neumann boundary conditions on $\bd_1 M_j$ and vanishes on $\bd_2 M_j$. As $M_j$ is strictly included in $M_{j+1}$, by the unique continuation principle of Aronszajn (see \cite{Aro57}), the sequence $\left(\lambda_0(M_j)\right)_{j\in\Bb N}$ is strictly decreasing and $$\forall j\in\Bb N,\ \lambda_0(M_j)>\lambda_0(M).$$

Let $p^j(x,y,t) = p_{M_j}(x,y,t)$ be the heat kernel of $M_j$ associated to our problem, that is the fundamental solution of the equation
$$\Delta f = -\frac{\bd f}{\bd t}$$
with these mixed boundary conditions (recall our Laplacian is positive). For any $x,y\in M_j$ and $t>0$,
\beq \label{eq:NoyChal} p^j(x,y,t) = \sum_k e^{-\lambda_k^j t}\phi_k^j(x)\phi_k^j(y),\eeq
where $\phi_k^j$ is the eigenfunction of the Laplacian with Neumann boundary conditions on $\bd^1M_j$ and Dirichlet boundary conditions on $\bd^2M_j$ associated to the eigenvalue $\lambda_k^j$. 
\pgh

The following statements would need, to be proved, a long development of the \emph{diffusion theory} associated to an elliptic operator, which go through the theory of \emph{stochastic calculus}. Therefore, we will only present definitions, results and key ideas of the proof. We lead our reader to \cite{Mall75} and \cite{Tay2}, Chapter 11 for the construction of the Brownian motion from the heat kernel. Probabilstic basics for this proof can be found for example in \cite{Bas95}, chapter I, and the details of our proof from stochastic integrals are in \cite{Bas98} in the case of open domains in $\Bb R^d$. The justification of their adaptation for riemannian manifolds is for example in \cite{Em89}. 

\begin{Defi}
Let $j>0$ and $M_j\subset$ one of the domain we defined above. We note $\Omega$ the set of all continuous path from $\bar{\Bb R}_+$ into $M$ and $\tau_j : \Omega\ra \bar{\Bb R}_+$ defined for any $\omega\in\Omega$ by
\beq \label{eq:Tau} \tau_j(\omega) = \inf\left\{t>0 : \omega(t)\in\bd^2 M_j\right\}.\eeq
Let $\Omega^j$ be the set of path from $\bar{\Bb R}_+$ into $\bar M_j$ such that 
$$\forall t \geq\tau(\omega), \omega(t)=\omega(\tau(\omega))\in \bd^2 M_j$$
and
$$\Omega^j_x = \left\{\omega\in \Omega^j :\omega(0) = x\right\}.$$
We call $\Omega^j$ the set of \ind{trajectories} in $M_j$.
\end{Defi}

We call a \ind{cylinder} in $\Omega^j_x$ a set of the following form :
$$A = \left\{\omega\in \Omega^j_x : (\omega(t_1),\ldots,\omega(t_k))\in B\right\},$$
with $k\in \Bb N$, $B\subset (M_j)^k$ is a borelian, and the $t_j$ are real number such that $0\leq t_1<t_2...<t_k.$ For any cylinder $A$, we set
$$\Bb P^j_x(A) = \int_Bp^j(x,y_1,t_1)p^j(y_1,y_2,t_2-t_1)...p^j(y_{k-1},y_k,t_k-t_{k-1})dV(y_1)...dV(y_k),$$
where $dV$ is the canonical riemannian volume measure. One can show from the semi-group property of the heat kernel that $\Bb P^j_x$ extends to a unique \emph{probability measure} on the $\sigma$-algebra of $\Omega^j_x$ generated by its cylinders (see \cite{Tay2}). $\Omega^j$ is the \ind{Wiener space} on $M_j$, and $\Bb P_x^j$ the \ind{Wiener probability measure} in $x$. 

Let us consider the random process $(X_t)_{t\geq0}$ on $(\Omega_x^j,\Bb P^j_x)$ defined for any $\omega\in\Omega_x^j$ by
$$X_t(\omega) = \omega(t).$$
By definition of $\Bb P_x^j$, for any borelian $B$ of $M$ we have
$$\Bb P_x^j(X_{t+s}\in B|X_t = z) = \int_Bp^j(z,y,s)dV(y) = \Bb P_z^j(\{X_s\in B\})\ :$$
$X_t$ is a Markov process with law $p^j$.

The random process $(X_t)_{t\geq 0}$ with law $\Bb P_x^j$ is called the \ind{brownion motion} on $M_j$ (with reflexion on $\bd^1 M_j$, which will be assumed from now on). This name comes from the fact that, in a statistical physic model as the perfect gaz, $p(x,y,t)dV(y)$ is the density of probability for a particule that was in $x$ at $t=0$ to be in $y$ at time $t$.

Let $f : M^j\ra \Bb R$ be a $\Cl C^2$ function. For any vector $Y,Z\in T_xM$, we write now $g_x(Y,Z) = Y.Z$. The Itô Formula for the Brownian Motion $(X_t)_{t\geq 0}$ (see \cite{Bas95} p 49, \cite{Em89} p 34) is now :

\beq \label{eq:Ito1} f(X_t) = f(X_0) + \int_0^t \nabla f(X_s).dX_s-\int_0^t\Delta f(X_s)ds.\eeq

The las term of our formula is multiplied by $-2$ from the one in \cite{Bas95} p 49 : this comes from the fact that our sign convention for the Laplacian is opposed to the convention chosen by Bass, and that the usual Brownian motion considered by probabilists has elementary transition probability given by the fundamental solution of
$$-\frac{1}{2}\Delta f = \frac{\bd f}{\bd t},$$
whereas, as geometers, we do not keep this $\frac{1}{2}$ in our definition (see for instance \cite{Bas98} p53).

By \cite{Bas98} p 33, as $X_t$ is a brownian with normal reflexion on $\bd^1 M_j$, we can write
$$dX_t = dW_t + \nu(X_t)dL_t,$$
where $W_t$ is a brownian without reflexion in $M_j$, $\nu(X_t)$ is the inward norml to $\bd^1M_j$ in $X_t$ when $X_t\in\bd^1M_j$, and $0$ elsewhere, and $L_t$ is the \ind{local time} on $\bd^1M_j$. This local time is a positive increasing process with bounded variations, strictly increasing if and only if $X_t\in\bd^1M$. It is defined by
$$L_t = \lim_{\epsilon\ra 0} \frac{1}{\epsilon}\int_0^t \textbf{1}_{d(X_s,\bd^1M_j)\leq\epsilon}ds,$$
where $d(X_s,\bd^1M_j)$ is the distance from $X_s$ to $\bd^1M_j$ induced by the metric on $M$. The formula (\ref{eq:Ito1}) becomes then

\beq \label{eq:Ito2} f(X_t) = f(X_0)+\int_0^t \nabla f(X_s).dW_s+\int_0^t\nabla f(X_t).\nu(X_t)dL_t-\int_0^t\Delta f(X_s)ds.\eeq

\subsection{A lower bound from superharmonic functions}

A smooth function is said to be $\lambda$-surharmonic if and only if it satisfies $(\Delta f\geq \lambda f$.
A first half of Theorem \ref{th:SullNeum} will be implied by the following proposition :

\begin{Prop}\label{prop:SurharmPos}
Let $\lambda\in\Bb R$ such that there exists a positive $\lambda$-superharmonic function on $M$ with Neumann boundary conditions on $\bd M$. Then $\lambda_0(M)\geq \lambda$.
\end{Prop}
\begin{proof}

Let $j\in\Bb N$, for any $\Cl C^2$ positive function $f$ and real $\lambda$, we consider the random process $(Y_t)_{t\geq0}$ on $\Omega_x^j$ with real values defined by
$$Y_t = e^{\lambda \min(t,\tau_j)} f(X_t),$$
where $X_t$ is still the brownian motion with reflexion we defined in the previous section and $\tau_j$ is the first time $X_t$ reaches $\bd_2 M_j$. Itô Formula implies now
\beq \label{eq:Ito3} e^{\lambda \min(t,\tau_j)} f(X_t) = f(X_0) +\int_0^t e^{\lambda \min(s,\tau_j)} \nabla f(X_s).dX_s+ \int_0^t \lambda e^{\lambda \min(s,\tau_j)} f(X_s)ds -\int_0^te^{\lambda \min(s,\tau_j)}\Delta f(X_s)ds.\eeq

Assume now $\Delta f\geq \lambda f$ and $f$ satisfies Neumann boundary conditions on $\bd_1 M_j$. Note $\tau_j$ satisfies
$$\Bb P_x^j(\tau_j>t) = \int_M p^j(x,y,t)dV(y).$$
Since
$$p^j(x,y,t) = \sum_k e^{-\lambda_k^j t}\phi_k^j(x)\phi_k^j(y),$$ we have
$$\lim_{t\ra\infty} e^{\lambda_0^j t}p^j(x,y,t) = \phi_0^j(x)\phi_0^j(y).$$
As $\phi_0>0$ on $\inter{M}_j$,
\beq \label{eq:TFinPS} \Bb P^j_x(\{\tau_j> t\}) \sim C e^{-\lambda_0^j t}. \eeq
Therefore, for any $\lambda<\lambda_0^j$,   $e^{\lambda_0}\tau_j$ is almost surely finite and $Y_t$ is summable on $\Omega_x^j$. Hence, for any $t>0$, integrating the formula (\ref{eq:Ito3}) on $\Omega^j_x$, we get
$$f(x) \geq \Bb E_x^j(e^{\lambda t}f(X_t)) = \int_{M_j} e^{\lambda t}f(y)p^j(x,y,t)dV(y) + \int_{\Omega^j_x} e^{\lambda\tau(\omega)}f(\omega(\tau(\omega)))\textbf{1}_{\{\tau\leq t\}}d\Bb P_x^j(\omega).$$
Hence, for any $t>0$,
$$f(x)\geq \int_{M_j} e^{\lambda t}f(y)p^j(x,y,t)dV(y).$$
As $f(x)$ is finite and $$\lim_{t\ra\infty} e^{\lambda_0^j t}p^j(x,y,t) = \phi_0^j(x)\phi_0^j(y),$$ with $\phi_0^j>0$ on $\inter{M_j}$, then we must have $\lambda\leq\lambda_0^j$. As this is true for any $j\in\Bb N$ and $\lambda_0(M) = \sup_j \lambda_0(M_j)$ this concludes the proof of our proposition.
\end{proof}

To prove Theorem \ref{th:SpecDom} in Section 2, we only needed this half of Theorem \ref{th:SullNeum}. We now sketch a proof of the other half for completeness.

\subsection{Existence of a $\lambda_0$-harmonic function}

To finish the proof of Theorem \ref{th:SullNeum}, it is enough to show there exists a $\lambda_0$-harmonic positive function on $M$.

The random variable $\tau_j : \Omega\ra \bar{\Bb R}_+$ defined in (\ref{eq:Tau}) is a \ind{stopping time} (see \cite{Bas95} p 13) which is almost surely finite by (\ref{eq:TFinPS}). Now, by Doob Stopping Theorem (see \cite{Bas95} p 29), integrating Formula \ref{eq:Ito3} gives :

\begin{Theo}\label{th:LHarm2}
For any $\lambda<\lambda_0^j$ and any $\lambda$-harmonic positive function $f$ on $M^j$ with Neumann condition on $\bd^1 M^j$,
$$f(x) = \Bb E_x^j(e^{\lambda \tau}f(X_\tau)) = \int_{\Omega_x^j}e^{\lambda\tau(\omega)}f(\omega(\tau(\omega)))d\Bb P_x^j(\omega) = \int_{\bd M_j}f(\xi)d\mu^\lambda_{j,x}(\xi),$$

where $\mu_{j,x}^\lambda$ is the (finite, non normalized) measure defined for any borelian $B$ of $\bd_2 M_j$ by :
$$\mu_{j,x}^\lambda(B) = \Bb E_x^j(e^{\lambda\tau} \textbf{1}_{\{X_\tau\in B\}})= \int_{\Omega_x^j} e^{\lambda \tau(\omega)} \textbf{1}_{\{\omega(\tau(\omega))\in B\}} d\Bb P_x^j(\omega).$$
\end{Theo}

Now, for any $\lambda<\lambda_0(M_j)$ this formula can be used to extend any function on $\bd_2 M_j$ to a $\lambda$-harmonic function on $M_j$:

\begin{theo}\label{th:ExLHarm}
Let $f : \bd_2 M_j\ra\Bb R^+$ be a borelian non-negative function on $\bd_2 M_j$, non identically 0. Then for any $\lambda\leq\lambda_0(M)<\lambda_0(M_j)$, the function $\tilde f : M_j\ra \Bb R^*_+$ defined for any $x\in M_j$ by
$$\tilde f(x) = \int_{\bd M_j}f(\xi)d\mu^\lambda_{j,x}(\xi)$$
is $\lambda$-harmonic, positive on $M_j$, satisfies Neumann conditions on $\bd^1 M_j$ and is continuous on $\bar M_j$.
\end{theo}

When $\lambda=0$, $\mu_{j,x}$ is called the \ind{Poisson measure} or \ind{harmonic measure} on $\bd^2 M_j$ (with reflexion on $\bd^1 M_j$) starting from $x$. This last theorem can be proven from the properties of the measures $\mu^\lambda_{j,x}$ which come from the heat kernel they were built with. We will not do it here, the reader can refer to Section 5 of Chapter 11 of \cite{Tay2}.
\pgh
As for any $j\in\Bb N, \lambda_0 = \lambda_0(M)<\lambda_0(M_j)$, by this theorem we construct a sequence of $\lambda_0$-harmonic functions $(f_j)_{j\in\Bb N}$ with Neumann boundary condition on $\bd M$ such that for any $j\geq 0$, $f_j$ is positive on $M_j$. Let $x\in M_0$, we can assume for all $j\in\Bb N, f_j(x) = 1$. Then, using a \ind{Harnack Principle}, one can show that these $f_i$ converge uniformly on any compact subset of $M$ to a smooth $\lambda_0$-harmonic function $f_\infty$, positive on $M$ and satisfying Neumann boundary conditions on $\bd M$. The proof of this last argument is verbatim in \cite{Sul87}, p 336.

\bibliographystyle{alpha}
\def\cprime{$'$}

\end{document}